\documentclass[10pt]{article}
\usepackage[dvips]{color}
\usepackage [autostyle, english = american]{csquotes}
\usepackage[
pdftitle={A fast, high-order scheme for volume potential evaluation},
pdfcreator={pdflatex},
pdfsubject={preprint},
hyperindex = {true},
colorlinks = {true},
linkcolor = {blue},
citecolor = {blue}
]{hyperref}

\usepackage[backend=biber,
  style=numeric,
  sorting=none,
  maxbibnames=99,
  uniquename=false,
  giveninits=true,
  uniquelist=false]{biblatex}
\bibliography{bibliography}
\renewbibmacro{in:}{%
  \ifentrytype{article}
    {}
    {\bibstring{in}%
     \printunit{\intitlepunct}}}

\usepackage{tikz}
\usetikzlibrary{shapes,decorations}
\usepackage[english]{babel}
\usepackage{indentfirst}
\usepackage{amsmath}
\usepackage{algorithm,algpseudocode}
\usepackage{wrapfig}
\usepackage{bm}
\usepackage{amssymb}
\usepackage{caption}
\usepackage[shortlabels]{enumitem}
\usepackage{booktabs}
\usepackage{authblk}
\usepackage{graphicx}
\usepackage{amsthm}
\usepackage{comment}
\usepackage{cleveref}
\newtheorem{theorem}{Theorem}
\newtheorem{lemma}{Lemma}
\newtheorem{corollary}{Corollary}
\theoremstyle{remark}
\newtheorem{remark}[theorem]{Remark}

\usepackage[title]{appendix}
\usepackage{float}
\usepackage[lofdepth,lotdepth]{subfig}

\DeclareMathAlphabet{\mathpzc}{OT1}{pzc}{m}{it}
\usepackage{xargs}
\usepackage{matlab-prettifier}
\usepackage[nodisplayskipstretch]{setspace}
\numberwithin{equation}{section}
\usepackage{wrapfig}
\usepackage{titling}
\graphicspath{{fig/}{../fig/}}
\newcommand{\vv}[1]{\boldsymbol{#1}}
\renewcommand{\d}{{\mathrm{d}}}

\allowdisplaybreaks

\DeclareMathOperator*{\argmin}{argmin}

\DeclareRobustCommand{\rchi}{{\mathpalette\irchi\relax}}
\newcommand{\irchi}[2]{\raisebox{\depth}{$#1\chi$}} 

\title{A fast, high-order scheme for evaluating volume potentials on complex 2D geometries via area-to-line integral conversion and domain mappings}

\author{Thomas G. Anderson\thanks{Department of Mathematics, University of Michigan},\, Hai Zhu\textsuperscript{*}\hspace{-.4em},\, and Shravan Veerapaneni\textsuperscript{*}}
\date{}

\begin{document}
\maketitle

\vspace{-2em}

\begin{abstract}
    While potential theoretic techniques have received significant interest and found broad success in the solution of linear partial differential equations (PDEs) in mathematical physics, limited adoption is reported in the case of nonlinear and/or inhomogeneous problems (i.e.\ with distributed volumetric sources) owing to outstanding challenges in producing a particular solution on complex domains while simultaneously respecting the competing ideals of allowing complete geometric flexibility, enabling source adaptivity, and achieving optimal computational complexity. This article presents a new high-order accurate algorithm for finding a particular solution to the PDE by means of a convolution of the volumetric source function with the Green's function in complex geometries. Utilizing volumetric domain decomposition, the integral is computed over a union of regular boxes (lending the scheme compatibility with adaptive box codes) and triangular regions (which may be potentially curved near boundaries). Singular and near-singular quadrature is handled by converting integrals on volumetric regions to line integrals bounding a reference volume cell using cell mappings and elements of the Poincar\'e lemma, followed by leveraging existing one-dimensional near-singular and singular quadratures appropriate to the singular nature of the kernel. The scheme achieves compatibility with fast multipole methods (FMMs) and thereby optimal asymptotic complexity by coupling global rules for target-independent quadrature of smooth functions to local target-dependent singular quadrature corrections, and it relies on orthogonal polynomial systems on each cell for well-conditioned, high-order and efficient (with respect to number of required volume function evaluations) approximation of arbitrary volumetric sources. Our domain discretization scheme is naturally compatible with standard meshing software such as \texttt{Gmsh}, which are employed to discretize a narrow region surrounding the domain boundaries. We present 8th-order accurate results, demonstrate the success of the method with examples showing up to 12-digit accuracy on complex geometries, and, for static geometries, our numerical examples show well over $99\%$ of evaluation time of the particular solution is spent in the FMM step.
\end{abstract}

\section{Introduction}
\label{sec:intro}

This article describes a fast, high-order accurate numerical scheme for evaluating the {\em volume potential}, also known as the {\em Newton potential}, given by
\begin{equation}\label{eq:volPot}
    \mathcal{V}[f]\left(\vv{r}_0\right) = \int_{\Omega} G(\vv{r}, \vv{r}_0)\, f(\vv{r})\, \mathrm{d}A(\vv{r}), \quad \vv{r}_0 \in \Omega,
\end{equation}
where $\Omega$ is an irregular two-dimensional domain, $f$ is a given source density function and $G$ is typically the free-space Green's function for an underlying linear, constant-coefficient elliptic partial differential equation (PDE) operator, such as for the Laplace, Stokes or (modified) Helmholtz equations. Domain convolutions of this form are a key ingredient when solving inhomogeneous PDEs {\em via} potential theory \cite{hsiao2008boundary}. Historically, the community has tackled a variety of fundamental challenges in the use of potential theoretic methods for addressing such problems.

While integral equation methods offer a reduction in dimensionality in the associated homogeneous problem, they lead to dense operators upon discretization, whose efficient application has formed a significant body of work over the recent decades (such as the FMM~\cite{greengard1987fast}).
With regards to discretization, as the kernels which arise in boundary integral formulations for homogeneous PDE boundary value problems are singular, the integral equations that result involve singular and nearly-singular integrals. Singular quadrature has formed the basis of significant inquiry with many successful schemes proposed; see~\cite{Hao:14} for a recent review of the situation in two dimensions (note that in three dimensions the landscape is more challenging). Furthermore, {\em nearly singular} integrals arise from the nature of these kernels in the presence of arbitrary target points close to $\partial \Omega$, as generally occur in complex geometries e.g.\ when boundary components are close to touching. These challenges each have a counterpart in the volume potential problem. Recently, a concerted effort in the field produced several robust strategies for high-accuracy evaluation of nearly singular integrals in two dimensions \cite{helsing2008evaluation, klockner2013quadrature, barnett2015spectrally, rahimian2018ubiquitous, af2018adaptive, wu2020solution,Perez:19,Faria:21}. For example, the method of \cite{barnett2015spectrally}, which we employ for the homogeneous solver component of this work, achieves spectral accuracy in evaluating layer potentials close to smooth boundaries, requiring only a modest number of discretization nodes per domain inclusion for full accuracy in double precision. These advances motivate us to revisit the volume potential evaluation problem as a prominent remaining task.

Most fundamentally for inhomogeneous problems, the volume discretization needs to be adapted to the irregular domain.  We first mention that for the case of $\Omega$ a unit box or a hierarchical finite union of scaled and translated boxes, there exist fast, adaptive, highly-accurate and computationally scalable algorithms~\cite{ethridge2001new, langston2011free, malhotra2014volume, malhotra2015pvfmm}. In a similar vein, reference~\cite{greengard1996direct} solves local elliptic problems on each box using the action of the operator on orthogonal polynomials and obtains high-order accuracy with very limited computational expense per degree of freedom. Taken together, these ``box code'' methods have proven highly effective because they inherently exploit the translational invariance of the Green function of the PDE, allowing repeated use of what are essentially lookup tables for singular and near-singular evaluation points within a vicinity of a given box. However, achieving the same level of accuracy and efficiency when $\Omega$ is an arbitrary complex geometry has been a sustained challenge. Na\"ively, for a box code in the presence of complex geometry some boxes inevitably are `cut' by the boundary and evaluating the contribution to the volume potential~\eqref{eq:volPot} from such irregular `cut cells' is the primary challenge to be met (such issues arise in other contexts utilizing regular grids in the presence of embedded boundaries see e.g.~\cite{Giuliani:21}). Attempts have been made to address irregular geometries by using (generally low-order accurate) extrapolation or local extension for the function $f$ over the resulting cut cells in combination with extensive adaptivity near to domain boundaries to achieve desired tolerances~\cite{langston:12thesis}. As a result, while the resulting methods effectively make use of fast algorithms to reduce the computational burden, the required number of degrees of freedom appears to be significant.

In seeking to overcome this issue, most existing works avoid direct evaluation of \eqref{eq:volPot} in an irregular domain; instead, one class of algorithms employ (local) volumetric PDE solvers---finite difference~\cite{Mayo:84,Mayo:92,Mayo:95,Rapaka:20} and finite element~\cite{biros2004fast} methods---in an embedded regular domain and layer potentials for enforcing correct boundary data. In particular, the embedded boundary integral approach of~\cite{biros2004fast} solves the inhomogeneous PDE problem on a rectangular domain that embeds $\Omega$; local low-order extrapolation for bulk forces and use of jump relations result in a second-order accurate method, which appears quite efficient.
The method of reference~\cite{Rapaka:20}, meanwhile, couples to geometric multigrid solvers for the bulk and is also restricted to the Poisson equation; this work in fact includes a timing comparison to the box code method of reference~\cite{askham2017adaptive} (discussed in the next paragraph), which shows superior speed per degree of freedom. We note, however, that only first-order convergence is demonstrated in~\cite{Rapaka:20}, with no claim made that the accuracy in the comparison example matches that of the box code with continuous extension, and we further observe that while the method implicitly requires extension of the source density outside the domain, this point is unaddressed (i.e.\ known source functions are chosen that are continuous across the boundary and their values outside used).

Smooth extension or continuation methods for the volumetric source form another broad class of methods, in which the irregular domain $\Omega$, again, is embedded in a box $B$, the function $f$ is extended onto $B  \setminus \Omega$, and standard box codes are applied on the whole of $B$ (inevitably incurring some performance penalty due to the increased degrees of freedom over the enlarged domain).  It is essential in these methods that the extended function be smooth for high-order accuracy. In~\cite{askham2017adaptive}, a harmonic extension of $f$ is performed using a boundary integral approach and is coupled to an adaptive box code.
While the approach generalizes beyond the Poisson equation solver described in this work, its main limitation is that constructing higher-order smooth approximations requires higher-order derivatives of the source function $f$ along the domain boundary. More significantly, it requires solution of a separate harmonic equation in order to generate the continuously extended source function, with higher-order continuity for higher-order convergence requiring solution of progressively more complex and costly high-order harmonic equations. Indeed the use of adaptivity in~\cite{askham2017adaptive} functions at least in part as a means to achieve high accuracy in the presence of limited smoothness in the extension, as experienced previously in e.g.~\cite{ethridge2001new,langston2011free,langston:12thesis}.

Other approaches to extension have traded off adaptivity for high-order accuracy and compatibility with FFTs on uniform grids; one recently introduced technique is the two-dimensional Fourier continuation method~\cite{bruno2020two} which as well as being a general-purpose numerical tool has recently been demonstrated to give accurate smooth periodic extensions suitable for use in elliptic solvers, at least in some simple geometrical contexts. In~\cite{fryklund2018partition}, a function extension scheme based on the use of radial basis functions (RBFs) is proposed, termed as the {\em partition of unity extension} (PUX) method. PUX lays down a set of disks along the domain boundary, each covering uniform points both inside and outside the domain, and solves overdetermined linear systems on each disk to generate smooth extensions to Cartesian grid points within these disks laying outside the domain, which are then smoothly mollified to zero away from the boundary. Although PUX has shown success on a variety of examples~\cite{af2020fast,fryklund2020integral}, several hurdles remain for a scalable implementation including the number of parameters and the heuristic nature by which they are determined (e.g., the partition radius and shape parameters of the RBFs), poor conditioning of RBF methods and a reliance on uniform grids. 

\begin{figure}[!b]
    \centering
    \includegraphics[height=0.24\textwidth]{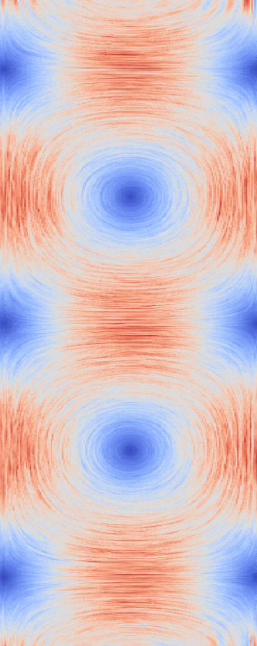}
    \includegraphics[height=0.24\textwidth]{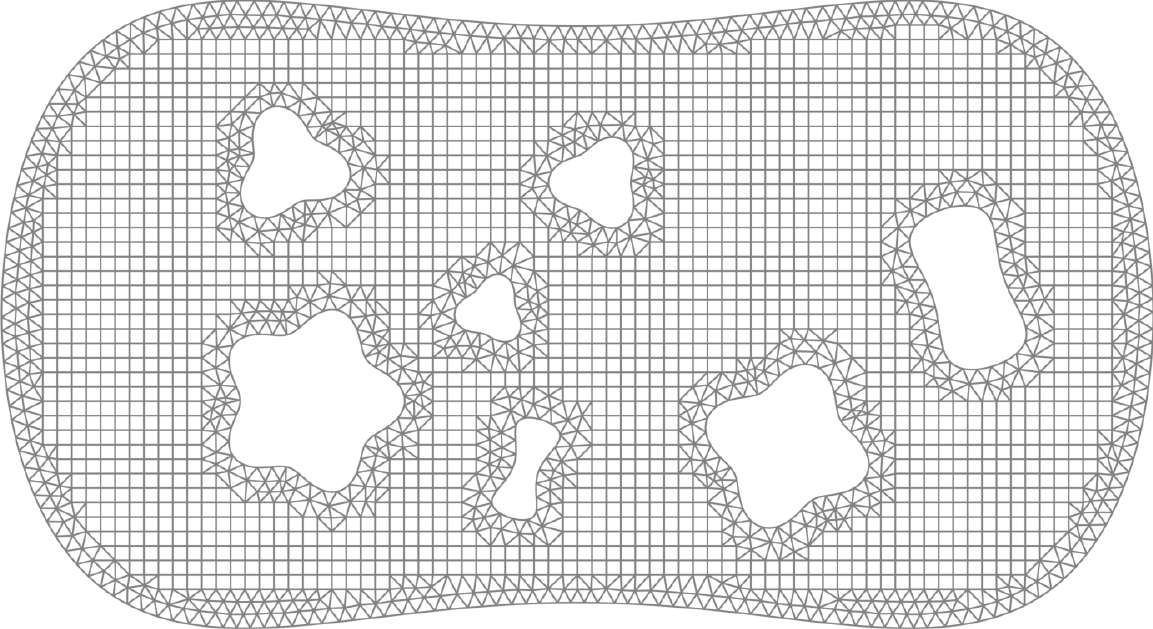}
    \includegraphics[height=0.24\textwidth]{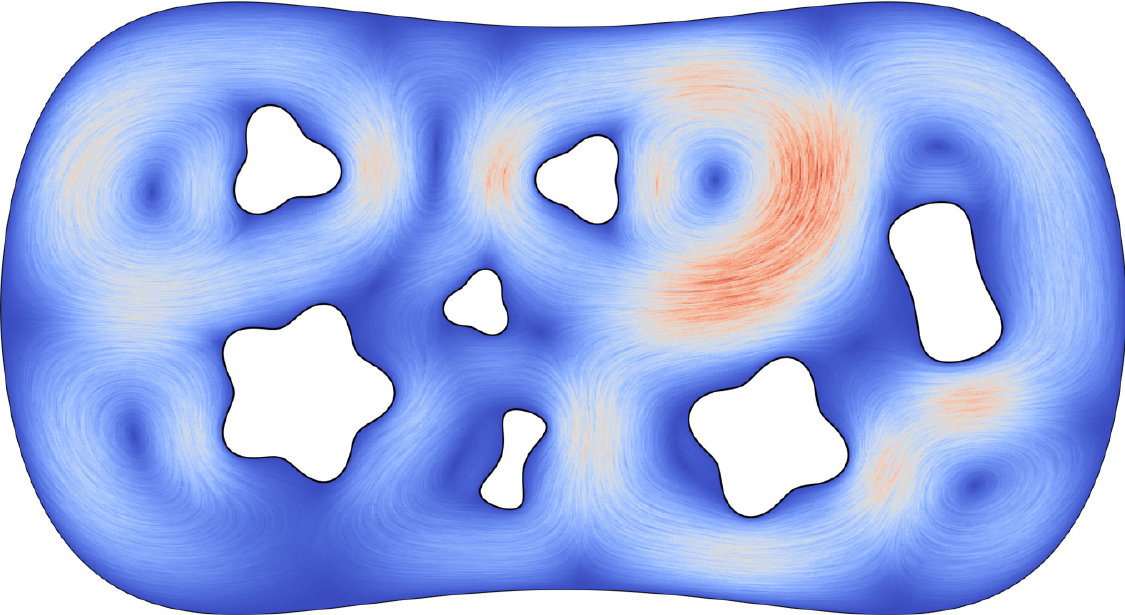}
    \caption{{\em\small Solution of the inhomogenous Stokes equations~\cite{hsiao2008boundary} on an irregular domain with no-slip boundaries. Left: contour plot of the imposed forcing function, $\vv{f}(x,y) = \left(-\sin x \cos y,\, \cos x\sin y\right)^T$. Center: volume mesh employed in our Newton potential \eqref{eq:volPot} evaluation scheme comprised of a regular grid in the bulk and a narrow triangulated region around the domain boundaries. Right: streamlines of the solution with color indicating the magnitude of the velocity field.
    }}
    \label{fig:Stokes}
\end{figure}

In this article, we take a more direct approach to solving inhomogeneous linear PDEs, by computing explicitly the volume integral~\eqref{eq:volPot} by means of numerical quadrature over complex geometry (an example is shown in Figure~\ref{fig:Stokes}). This mathematically obvious avenue has in the past been considered computationally impractical, with recent contributions~\cite{askham2017adaptive,af2020fast,fryklund2020integral} noting that such an approach is to be avoided because quadratures for irregular cells are seemingly challenging and less amenable to fast algorithms (reference~\cite{Steinbach:10} does evaluate~\eqref{eq:volPot}, in a method coupled to the FMM). Indeed, multidimensional quadratures over general regions are much less developed (though see~\cite{Vioreanu:14,XiaoGimbutas:10,Xiao:10}) though there exist ad-hoc quadratures over a variety of specific regions~\cite{Dunavant:85,Cowper:73,Cools:93,LynessCools:94,Stroud:71}; see~\cite{Cools:97} for an extensive review. This relative lack of development holds especially so for singular integrands; in this more specialized context, in addition to the older work~\cite{LynessCools:94} on triangles we find the work~\cite{Strain:95,Aguilar:02} for locally-corrected quadratures on rectangles and work~\cite{Bremer:12, Bremer:13} for discretizations of boundary integral operators on mapped triangles that form surfaces in three dimensions.

To address the challenges of multi-dimensional quadrature over potentially arbitrary regions, a variety of works have converted a volumetric integral of interest \emph{with a smooth integrand} to integration along domain boundaries (potentially after some degree of domain-decomposition), after which numerical quadrature rules are used for various resulting one-dimensional integrals. Much of this work is focused on cubature for polyhedra, where we note contributions for this purpose~\cite{Sommariva:07,Sudhakar:14} based on theorems of vector calculus. One approach~\cite{Gunderman:21} treats volumes bounded by rational parametric curves using Green's theorem, while others~\cite{Muller:13, Saye:15} produce volume quadratures over implicitly-defined surfaces and volumes by a recursive dimensional-reductional algorithm and ultimately result in one-dimensional integrals that can be treated with Gaussian quadrature.

We also note a variety of volumetric-via-surface-integral quadratures that, like ours, use ideas from classical proofs of the Poincar\'e lemma (which, in brief, is applicable to domains which are star-shaped with respect to \emph{some} interior point---but not necessarily to all points---and leads to representations of volume integrals as surface integrals with an iterated integral whose associated physical integration points lay on rays from the boundary to the star-point). Recent work of~\cite{Barnett:21} results in integral representations similar in some respects to our own but only smooth integrands are considered there.
Other works such as~\cite{Fata:12} have focused on addressing the Newton potential problem specifically, but a target-centered coordinate system with the star-point as the origin results in quadrature points that lay outside the domain for regions that are non-convex (i.e.\ containing regions, and therefore potential evaluation points, with respect to which the domain is not star-shaped). This fact implicitly imposes upon such methods a reliance on function extension in order to obtain values of the source density $f$ at non-physical quadrature nodes (or, alternatively, moment matching)---such tasks (which, in view of previous discussion on function extension can be seen to be challenging) were not addressed in those works, instead relying on assumed-known analytic expressions outside of the domain. Related efforts to avoid limitations of this kind with schemes designed to rely on points only inside the computational domain are described in~\cite{Wang:17a, Wang:17b} (and include acceleration via FMM), but numerical results demonstrate first-order convergence (in the simplest cases) or no convergence at all. The dual reciprocity method (see~\cite{Partridge} for references to this body of literature and see~\cite{Gao:02,Gao:05} for some related works) is another method popular in the engineering literature that employs surface integrals and delivers low to modest accuracy, with similar challenges as mentioned above. Our proposed scheme, while still based in spirit on the Poincar\'e lemma and utilizing surface integrals, does not suffer from any of these difficulties (cf. \Cref{rem:starshaped_vs_convex}): all quadrature nodes lay inside the domain and the scheme leads to high-order accuracy.

Our work generalizes the volume potential scheme developed in the second author's thesis \cite{Zhu:21thesis} (also see~\cite{zhu2021high}), wherein, Poincar\'e's lemma and recursive product integration rules were employed for a few pairs of PDE kernels and approximating bases. On the other hand, the methods proposed here generate numerical quadratures for regions arising from domain decomposition of the integral~\eqref{eq:volPot}, which are tailored to kernel functions with a variety of singular behaviors (i.e.\ not restricted to PDE Green function kernels) and lead to efficient algorithms that are compatible with fast algorithms (such as FMMs). An important goal is to retain compatibility with adaptive box codes~\cite{ethridge2001new, malhotra2015pvfmm} due to their unparalleled speed and maturity---which we achieve by triangulating only a small boundary-fitted region---but we also seek high-order convergence and speed for repeated application of the volume potential with different source densities $f$. A secondary objective is the use of components (e.g.\ meshing, choice of basis, quadrature generation schemes) that will generalize naturally to three dimensional volume quadratures.

{\em \textbf{Synopsis.}}
The proposed method proceeds by meshing a thin region near to the boundary using a combination of curvilinear (with one side conforming to the boundary represented as parametric curves) and straight triangles---the number of triangular regions growing \emph{linearly} as the mesh is refined with the number of regular boxes growing quadratically (for the uniform gridding of the bulk that we consider here). Each cell (or element) is represented as a map from a reference element (either of a simplex or a box, with straight triangles corresponding to a simple affine map) and all physical target points are characterized by their position in reference space relative to cells for the purposes of singular and near-singular quadrature. Singular and near-singular quadrature on these irregular domains (and indeed even over boxes), in turn, is performed with high-order accuracy by use of a Poincar\'e lemma-type idea that re-writes the integrals (in reference space) on standard triangles (or boxes) over the boundary of that domain, followed by the use of one-dimensional quadratures that are adapted to the singular behavior of the kernel function (typically, the PDE Green function). The method possesses optimal asymptotic complexity (costing $\mathcal{O}(N)$ to obtain the Newton potential at, say, all $N$ quadrature points in a domain) as a result of use of point-FMMs for summation of target-independent quadratures for smooth integrands and via coupling to inexpensive local singular corrections that can be represented as a linear map of small size.

{\em Advantages.}\/ Our proposed methodology builds naturally on existing and growing bodies of work concerning each of {\em i)} one-dimensional singular and near-singular quadrature as well as {\em ii)} interpolation and quadrature of smooth functions over convex polyhedra by means of orthogonal polynomials. It is fully compatible with adaptive schemes for spatially-concentrated sources. Problem geometry is exactly captured with our approach which incorporates the local boundary information and integrates it with existing software for high-quality mesh-generation. High-order accuracy is easy to achieve and requires only a one-dimensional quadrature rule adequate for the singularity of the given kernel function and a known orthogonal polynomial system on the reference cell. The previous two points (meshing and orthogonal polynomial systems) provide direct means to generalize this work to three-dimensions. Finally, the methods are straightforward to integrate with existing acceleration techniques.

{\em Limitations.}\/ At the present moment our solver discretizes regular portions of the domain with uniformly sized boxes, which can be inefficient for functions with significant variation. This issue is easily addressed since our scheme is designed to be drop-in compatible with adaptive box codes for the bulk region that are highly effective in this context. A more serious limitation is the cost of generating local singular corrections; while the regular bulk region does not contribute to these costs (so that asymptotically the set-up costs pertain only to thin boundary regions), and we present techniques in \Cref{sec:sing_corr_koornwinder} to manage the remaining computational burden, computation of these quantities still forms the majority of the up-front cost of the method. Lastly, we restrict our attention in this work to static geometries only.

{\em \textbf{Outline.}} This article proceeds in \Cref{sec:background} with necessary preliminaries that present boundary integral equation formulations for some common linear PDEs, demonstrates the role that volume potentials play in their solution, and then presents in Sections~\ref{sec:scheme} and \Cref{sec:sing_corr_koornwinder} the proposed methodology: \Cref{sec:scheme} covers volumetric meshing, then the use of smooth quadratures and interpolation on the coordinate-mapped triangles that naturally arise, followed by a description of novel singular and near-singular quadrature techniques, while \Cref{sec:sing_corr_koornwinder}, in turn, focuses on implementation and efficiency considerations of the method, both for repeated evaluation of~\eqref{eq:volPot} and for efficient up-front generation of necessary corrections to the smooth quadratures (as well as some limited discussion of our approach to the routine task of quadratures on the bulk region). A variety of numerical examples are presented in \Cref{sec:results} demonstrating the properties of the method and a brief summary with concluding remarks is given in \Cref{sec:Conclusions}.

\section{Potential theory for inhomogeneous PDEs}\label{sec:background}
This article is primarily concerned with methods to obtain solutions $u: \Omega \to \mathbb{R}$ to the elliptic PDE boundary value problem
\begin{subequations}\label{eq:ellipt}
\begin{align}
    L u (\vv{r}) &= f(\vv{r}),\quad \vv{r} \in \Omega,\label{eq:ellipt_prob_pde}\\
    u(\vv{r}) &= g(\vv{r}), \quad \vv{r} \in \Gamma\label{eq:ellipt_prob_bc},
\end{align}
\end{subequations}
where $\Gamma = \partial \Omega$ denotes the (assumed piecewise-smooth) boundary of a complex geometry $\Omega \subset \mathbb{R}^2$. Here, it is assumed that $L$ is a linear operator with a known translation-invariant Green function $G(\vv{r}, \vv{r}_0) = G(\vv{r} - \vv{r}_0)$, as occurs e.g. for the Poisson, Helmholtz, modified Helmholtz, and Stokes equations.
A standard solution technique for the boundary value problem~\eqref{eq:ellipt} is to exploit linearity and seek a particular solution $u_P$ to the problem
\begin{equation}\label{eq:inhomog_PDE}
    L u_P (\vv{r}) = f(\vv{r}),\quad \vv{r} \in \Omega,
\end{equation}
with no care given to boundary conditions for $u_P$ on $\Gamma$, and then subsequently solve the augmented boundary value problem
\begin{subequations}\label{eq:ellipt_augment}
\begin{align}
    L u_H (\vv{r}) &= 0,\quad \vv{r} \in \Omega,\label{eq:ellipt_augment_prob_pde}\\
         u_H(\vv{r}) &= g(\vv{r}) - u_P(\vv{r}), \quad \vv{r} \in \Gamma\label{eq:ellipt_augment_prob_bc},
\end{align}
\end{subequations}
whereby the solution to~\eqref{eq:ellipt} is given by
\begin{equation}
    u(\vv{r}) = u_H(\vv{r}) + u_P(\vv{r}), \quad \vv{r} \in \Omega.
\end{equation}

Depending on the operator $L$ and the domain $\Omega$ a variety of boundary integral equation formulations may be appropriate, each requiring their own numerical analysis---such questions are not of concern here. It suffices to note that all require boundary values of a particular solution $u_P$ which can be obtained by evaluation of the volume potential~\eqref{eq:volPot}.
This article is thus concerned with the efficient and accurate evaluation of this integral for targets $\vv{r}_0 \in \Gamma$ (to provide boundary values of $u_P$ as data for \Cref{eq:ellipt_augment_prob_bc}) and for $\vv{r}_0 \in \Omega$ (for evaluation of the full solution $u$ at desired evaluation points $\vv{r} \in \Omega$). Nevertheless, we briefly outline for definiteness the integral equation formulations used in this article.

\subsection{Boundary integral formulation}
Consider the Poisson equation, for which $L = -\Delta$. We utilize standard representation formulas for the solution $u_H$ at a point $\vv{r}_0$ in the domain $\Omega$, expressed in terms of the double-layer potential
\begin{equation}\label{eq:representation_poisson}
    u_H(\vv{r}_0) = \mathcal{D}[\varphi](\vv{r}_0) := \int_\Gamma \frac{\partial G(\vv{r}, \vv{r}_0)}{\partial \vv{n}(\vv{r})} \varphi(\vv{r}) \,\d\sigma(\vv{r})
\end{equation}
induced by the boundary integral density $\varphi$, where $G(\vv{r}, \vv{r}_0)$ denotes the Green function of \Cref{eq:ellipt}, $G(\vv{r}, \vv{r}_0) = -\frac{1}{2\pi} \log\left(|\vv{r} - \vv{r}_0|\right)$.
As is well-known, using the representation formula~\eqref{eq:representation_poisson} and enforcing the boundary conditions~\eqref{eq:ellipt_augment_prob_bc} leads via the jump relations of the double layer potential~\cite{Kress} to the integral equation
\begin{equation}\label{eq:2ndkind_BIE_poisson}
    \begin{split}
        \left(\pm \frac{1}{2}I + D\right)[\varphi](\vv{r}) = g - u_P, \quad \vv{r} \in \Gamma_{\pm},
    \end{split}
\end{equation}
for a function $\varphi$ which must be satisfied in order for \Cref{eq:representation_poisson} to yield a solution to the boundary value problem~\eqref{eq:ellipt_augment} for $L = -\Delta$. Here $\Gamma_+$ (resp.\ $\Gamma_-$) denotes that section of the boundary $\partial \Omega$ with respect to which the domain lays exterior (interior), and we denote by $D$ the double-layer boundary integral operator
\begin{equation}
    D [\psi](\vv{r}_0) := \int_\Gamma \frac{\partial G (\vv{r}, \vv{r}_0)}{\partial \vv{n}(\vv{r})} \psi(\vv{r})\,\d\sigma(\vv{r}), \quad \vv{r}_0 \in \partial \Omega.
\end{equation}
The above formulation works as written for the modified Helmholtz equation ($L = -\Delta + \lambda^2$) as well by replacing the Green's function with   $G(\vv{r}, \vv{r}_0) := \frac{\lambda^2}{2\pi} K_0\left(\lambda |\vv{r} - \vv{r}_0|\right)$. We do not discuss here subtler points regarding treatment of nullspaces for these operators, see \cite{sifuentes2014randomized} for a treatise on this topic.

\section{A volume potential scheme for thin boundary-fitted regions}
\label{sec:scheme}
This section describes the basic elements of the proposed methodology for volume potential evaluation. It describes first an automatic method for the construction of a boundary-fitted region consisting of both curvilinear triangles (that are fitted to the boundary of inclusions defined by parametric curves) and straight triangles. Then, domain mappings are presented that map a reference cell to a triangular region and describe smooth quadrature and interpolation schemes for these regions (our approach to the routine task~\cite{ethridge2001new} of quadrature for regular source boxes are deferred and only briefly mentioned in \Cref{sec:boxes}). Finally, in \Cref{sec:sing_quad} we describe our approach to singular and near-singular quadrature based on expressing domain integrals over a reference cell in terms of line integrals on its boundary.

\subsection{Meshing}\label{sec:meshing}
This section describes a technique to generate a boundary-fitted mesh for a (possibly multiply-) connected volumetric region $\Omega$ with boundary formed by a collection of curves $\Gamma = \bigcup_{i=0}^{N_\Gamma} \Gamma_i$, with the method proceeding in an identical manner whether the domain is an unbounded (exterior) or bounded (interior) one.  The meshing algorithm we describe depends on parametrizations of the boundary curves, and to this end we introduce some useful notation.
We assume for each curve $\Gamma_i$ ($i = 1, \ldots, N_\Gamma$) that we have access to a global parametrization $\vv{\gamma}_i: [0, 2\pi] \to \Gamma_i$, $\vv{\gamma}_i = \vv{\gamma}_i(t)$, and denote by $\vv{n}_i$ the associated normal vector $\vv{n}_i = \vv{n}_i(t)$ directed into the domain $\Omega$. We also denote by $\ell_i(t)$ the arclength of the portion of the curve $\Gamma_i$ traced out by $\vv{\gamma}_i(\tau)$ for $0 \le \tau \le t$, and, abusing notation slightly, we call $\ell_i = \ell_i(2\pi)$ the total arclength of curve $\Gamma_i$. The algorithm fills the bulk of the domain away from the boundary with boxes in a uniform background mesh, assumed to be of size $h$.

The boundary-fitted mesh will be in the form of a tessellation involving $N_t$ triangular regions $\mathcal{T}_k$ and $N_b$ regular boxes $\mathcal{B}_k$,
\begin{equation}\label{eq:domain_tessellation}
    \Omega = \cup_{k=1}^K \mathcal{C}_k = \left(\cup_{k=1}^{N_t} \mathcal{T}_k\right) \bigcup \left(\cup_{k=1}^{N_b} \mathcal{B}_k\right),
\end{equation}
which together comprise $K = N_t + N_b$ cells, ordered first by triangles so that
\begin{equation}\label{eq:Ck_Tk_Bk}
    \mathcal{C}_k = \begin{cases}
    \mathcal{T}_k, \quad &1 \le k \le N_t,\\
    \mathcal{B}_{k-N_t}, \quad &N_t + 1 \le k \le K,
    \end{cases}
\end{equation}
where the box $\mathcal{B}_{k}$ has center $\vv{o}_k$.

Corresponding to the curve $\Gamma_i$ with arclength $\ell_i$ there will be a contribution of approximately $\ell_i/ h$ boundary-fitted (curved) triangles, reflecting general uniformity in the size of mesh cells that abut the boundary; other strategies incorporating considerations of local curvature or source adaptivity are also possible but beyond the scope of the present discussion (see also \Cref{rem:TriWild}).
For each of the $N_\Gamma$ curves, the algorithm proceeds in three steps to generate a boundary-fitted mesh: (1) Firstly, by identification of boundary `knot' points which segment each curve parametrization into $N_{\gamma_i}$ sections of approximately equal-arclength ($\Delta \ell \approx h$), then (2) Secondly, by generation of boundary-fitted mesh cells conforming to the knot points, and (3) Finally, by generation of ``buffer'' zones that connect the boundary-fitted mesh to the background ``bulk'' mesh.

\begin{figure}[!ht]
    \centering
    \includegraphics[width=\textwidth]{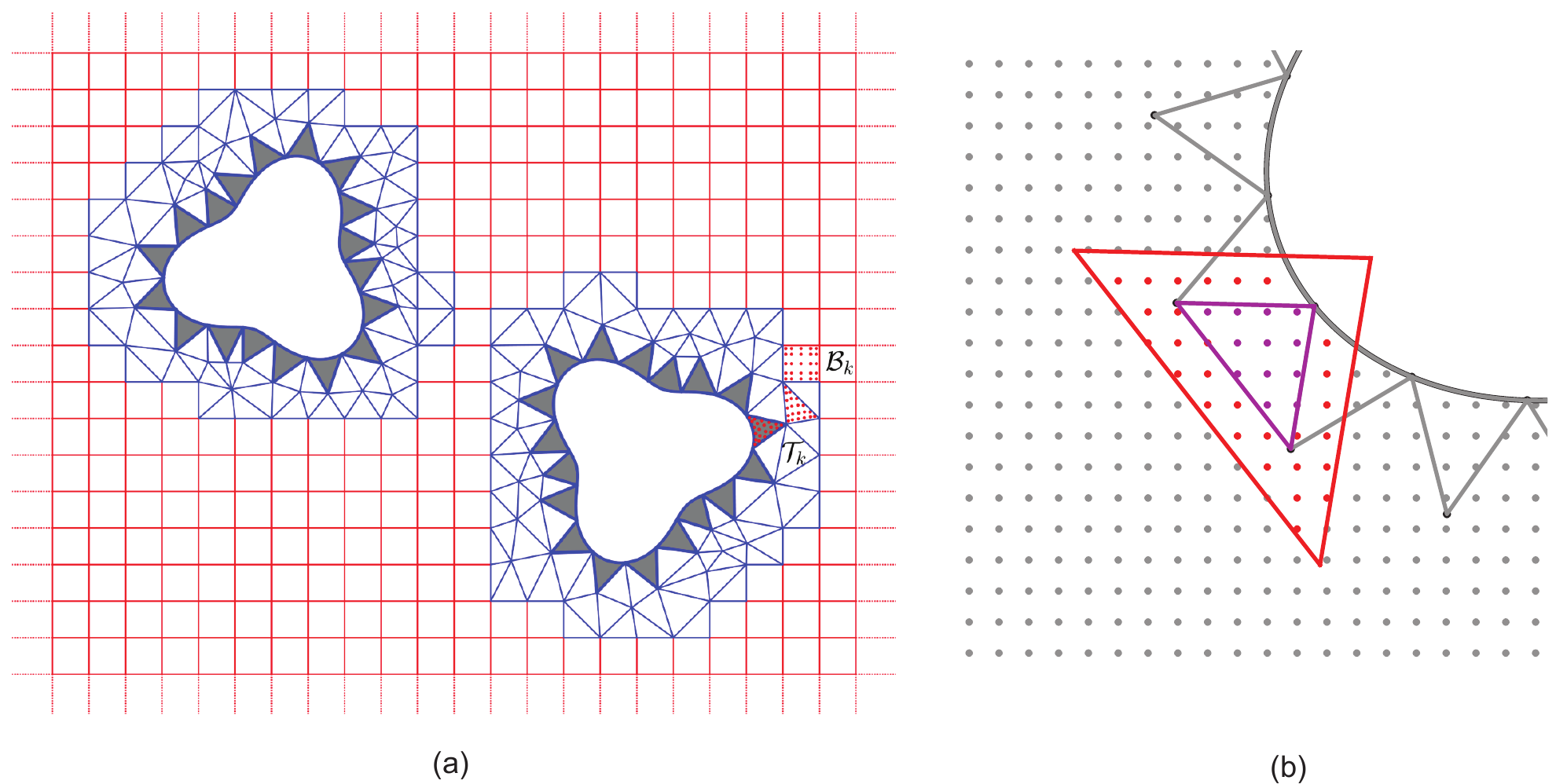}
    \caption{{\em\small (a) The volumetric discretization that arises from the present methodology applied to the exterior of two smooth bounded curves. Curved triangles that abut the curve are shaded while the buffer zone triangle cells are identified by blue edges. In three of the cells the interpolation nodes are denoted with red points. (b) The self-near-far separation rule for a straight triangular cell (with cell boundary plotted in purple and vertices in black) applied to uniformly-distributed volumetric targets. The region bounded by the red triangle contains all near-singular targets (marked red), the triangle itself bounds all singular targets (marked purple), while the remainder of targets (marked gray) are considered to be far quadrature to be treated with a smooth quadrature rule.}}\label{fig:mesh_demo}
\end{figure}
First, we describe a method for generating a sequence of parametric knot points $\kappa_j \in [0, 2\pi]$, $\lbrace \kappa_j \rbrace_{j=1}^{N_{\gamma_i}}$. Starting from $t=0$, knot points are progressively laid down so that each segment of the curve formed by consecutive knot points contains an arclength of approximately $\Delta \ell \approx h$; care is taken that all knot points are reasonably spaced (in particular, that knots $\kappa_j$ do not lay close to $t = 2\pi$---which may lead to close knot points in view of the periodicity of $\vv{\gamma}_i(t)$). Then, from each pair of knot points $\kappa_j$ and $\kappa_{j+1}$, $j = 1, \ldots, N_{\gamma_i} - 1,$ a curvilinear triangle is generated by connecting the two knot points by the boundary curve, and then further connecting to a third point in the volume generated by projecting a point a distance $h$ from the boundary in the normal direction $\vv{n}_i((\kappa_j + \kappa_{j+1})/2)$. The three points forming the new triangle thus are $\vv{\gamma}_i(\kappa_j)$, $\vv{\gamma}_i(\kappa_{j+1})$, and $\vv{\gamma}_i(m_j) + h\vv{n}_i(m_j)$, where $m_j = (\kappa_j + \kappa_{j+1})/2$ is the curve-parametric midpoint of the new triangle being formed. The curved cells that result from this procedure are shown in \Cref{fig:mesh_demo}, shaded grey.

Since the mesh cells formed in such a manner do not conform in any way to any uniform (or quadtree) background mesh, a watertight mesh is generated by means of a global `buffer' zone $\mathcal{F}$ that lays flush with these cells and serves to separate them from the background mesh. This region is defined by excluding any box from the background mesh laying within distance $h \delta$ of any boundary-fitting triangles (the selection $\delta = 0.8$ was made for the experiments in this article), the process repeating for each of the curves $\Gamma_i$, and yielding $N_\Gamma$ local buffer zones $B_i$ ($i = 1, \ldots, N_\Gamma$) which satisfy $\bigcup_i B_i = \mathcal{F}$. In the simple case that all such regions $B_i$ are pairwise disjoint, it is a routine meshing task to generate $N_\Gamma$ high-quality triangulations for each as each is bounded by two polygons; we use the \texttt{Gmsh}~\cite{Gmsh} software suite for this task. The process is repeated for each of the curves $\Gamma_i$, creating a total of $N_t$ (curved and straight) triangular mesh cells (again, for the simple case of pairwise disjoint $B_i$).

One possibility that must be addressed, then, is the situation that arises when `buffer' zones from two or more curves $\Gamma_i$ and $\Gamma_j$, $i \neq j$, are overlapping ($B_i \cap B_j \neq \emptyset$) as occurs naturally when $\Gamma_i$ and $\Gamma_j$ are not sufficiently well-separated. The algorithm handles this situation by utilizing a queue-like system that progresses iteratively through all remaining unmeshed curves to identify overlapping $B_i$, merging these, and triangulating the result. Selecting the first remaining curve $\Gamma_i$ and letting $\mathcal{F}_i = B_i$, for any curve $\Gamma_j$ ($j > i$) whose buffer zone $B_j$ satisfies $\mathcal{F}_i \cap B_j \neq \emptyset$ we let $\mathcal{F}_i = \mathcal{F}_i \cup B_j$ and remove $\Gamma_j$ from the queue of remaining curves. This process terminates when $j = N_\Gamma$ after which the index $i$ is increased to the next remaining curve and the process continues until no curves are remaining ($i = j = N_\Gamma$). The resulting sets $\mathcal{F}_i$ are pairwise disjoint and satisfy $\cup_i \mathcal{F}_i = \mathcal{F}$ (the union taken over all $i$ for which $\mathcal{F}_i$ is defined); each are triangulated individually, again with \texttt{Gmsh}. This procedure thus robustly handles the possibility of arbitrarily-many curves that lay close to each other.

A final depiction of the region surrounding one or more curves $\Gamma_i$ is given in \Cref{fig:mesh_demo} for well-separated curves (the buffer zone cells $\mathcal{F}$ are shaded blue), while examples of the buffer-zone merging process described in the previous paragraph can be seen in the mesh in \Cref{fig:Poisson_polydisp}. The set of regular cells is defined, finally, by the remaining part of the volume, per \Cref{eq:domain_tessellation}; in \Cref{fig:mesh_demo} these regular cells are colored red (note that regular cells may be enclosed by any given union-of-buffer-zones region $\mathcal{F}_i$ and are still handled by the box code).

A final element of the mesh construction process is a self-near-far separation rule, which classifies a given target $\vv{r}_0$ as a singular, near-singular, or smooth quadrature point with respect to each mesh cell. A target $\vv{r}_0 \not\in\mathcal{C}_k$ is considered a near-singular target of $\mathcal{C}_k$ if it lays within a polygon, with sides parallel to the physical cell (of course, this prescription is extended in an obvious manner to treat curvilinear triangles) and with sides laying a distance $D h$ (we use $D = 4/10$ for experiments in this article) from the boundary of $\mathcal{C}_k$. We denote by $\mathcal{C}^\mathrm{near}_k(\vv{r}_0)$ the index set of cells whose associated separation rule classify $\vv{r}_0$ as a near-singular target, with $\mathcal{C}^\mathrm{near}(\vv{r}_0)$ the union of all cells $\mathcal{C}_k$ with index in $\mathcal{C}^\mathrm{near}_k(\vv{r}_0)$. We denote by $\mathcal{C}^\mathrm{self}(\vv{r}_0)$ the cell to which $\vv{r}_0$ lays in, and define further $\mathcal{C}^\mathrm{far}(\vv{r}_0) = \Omega \setminus \left(\mathcal{C}^\mathrm{near}(\vv{r}_0) \cup \mathcal{C}^\mathrm{self}(\vv{r}_0) \right)$ as the remainder of $\Omega$ which is well-separated from $\vv{r}_0$ (see the depiction in Figure \ref{fig:mesh_demo}). Therefore, the volume potential can be decomposed as
\[
    \mathcal{V}[f](\vv{r}_0) = \left(\int_{\mathcal{C}^\mathrm{near}(\vv{r}_0)} + \int_{\mathcal{C}^\mathrm{self}(\vv{r}_0)} + \int_{\mathcal{C}^\mathrm{far}(\vv{r}_0)} \right) G(\vv{r}, \vv{r}_0) f(\vv{r})\,\mathrm{dA}(\vv{r})
\]
over near-singular, singular, and smooth regions of $\Omega$, respectively, for each target $\vv{r}_0$.

\begin{remark}\label{rem:param}
It is trivial to define for each triangular cell a boundary parametrization $Z_k: [0, 2\pi] \to \partial \mathcal{T}_k$, which could in principle be utilized by the same methodology outlined in \Cref{sec:sing_quad} to generate physical-space volumetric quadratures for targets located in or near $\mathcal{T}_k$. In practice, however, all quadrature will occur on the standard simplex $\widehat{\mathcal{T}}_0$ or the unit box $\widehat{\mathcal{B}}_0$ (see also \Cref{rem:starshaped_vs_convex}). It should be noted as well that the boundary of every box $\mathcal{B}_k$ in fact has the same parametrization, up to a translation, with similar implications for generating physical-space quadratures.
\end{remark}

\begin{remark}\label{rem:TriWild}
The prescription we give here for generating a high-quality volumetric mesh is certainly not the only method to obtain a compatible mesh for volume potential evaluation, but was effective for forming meshes for the examples in this article. Local curvature could be usefully incorporated into the algorithm described above in a blend of strategies based on arclength and changes in curvature, to achieve higher quality meshes. We mention as well the ``\texttt{TriWild}'' method of reference~\cite{Hu:19} for producing a valid volumetric mesh consisting of straight and curved triangles from a collection of curves whose union forms $\Gamma$, which has demonstrated success on a variety of real-world benchmark problems and is fully compatible with our curved triangle approach to the boundary-fitted region. (At the present moment implementations of the \texttt{TriWild} technique will result in only third-order approximations to the true boundary curve, motivating our exact approach to geometry representation, but it appears that this is not an essential limitation of the general \texttt{TriWild} technique; moreover, it appears straightforward to utilize slight perturbations of \texttt{TriWild} meshes and retain full geometric accuracy.)
\end{remark}

\subsection{Quadratures and interpolation for smooth integrands}\label{sec:sparse_func_approx}
In this section we describe first the overall approach to computing the volume potential $\mathcal{V}[f](\vv{r}_0)$ at arbitrary target points $\vv{r}_0$ via domain decomposition and then discuss quadrature of a potential that contains \emph{smooth} integrands (as well as the related source density interpolation problem) over a domain-decomposed cell $\mathcal{C}_k$. Efficient and accurate evaluation follows from a combination of (i) Standard fast summation technologies which, on the one hand allow for efficient computation of sums arising from target-independent quadratures and which are, on the other hand, necessarily inaccurate in the vicinity of a given target point $r_0$ for a singular kernel, as well as (ii) Accurate singular and near-singular \emph{target-dependent} quadrature \emph{corrections} applied locally in the vicinity of the given target point $r_0$. While the focus in this section is on point (i), the methodology does rely on a fixed orthogonal basis used to approximate the source function $f$ to high-order accuracy and for this reason we discuss also interpolation of smooth functions. Our approach is similar in spirit to that of reference~\cite{Greengard:21}.

We first sketch the unifying singular correction and functional approximation strategy. The Newton potential of a function $\rho$ at a given target point $\vv{r}_0 \in \Omega$ can be written as
\[
    \mathcal{V}[\rho](\vv{r}_0) = \sum_{k=1}^K \mathcal{V}_k[\rho](\vv{r}_0),\quad\mbox{where}\quad \mathcal{V}_k[\rho](\vv{r}_0) = \int_{\mathcal{C}_k} G(\vv{r}, \vv{r}_0) \rho(\vv{r}) \,\mathrm{dA}(\vv{r}).
\]

Quadrature will be performed using coordinate mappings $\vv{R}^k$ from a reference cell $\widehat{\mathcal{C}}_k$---either a unit square $\widehat{\mathcal{B}}_0$ or the unit simplex $\widehat{\mathcal{T}}_0$ (box regions are identical by translational invariance, so the use of mappings in this case is done merely for notational consistency). We will denote by
\begin{equation}\label{eq:Vkdef}
  \widehat{\mathcal{V}}_k[\rho](\vv{\zeta}_0) = \int_{\widehat{\mathcal{C}}_k} G\left(\vv{R}^k(\vv{\zeta}), \vv{R}^k(\vv{\zeta}_0)\right) J^k(\vv{\zeta}) \rho(\vv{\zeta})\,\mathrm{dA}(\vv{\zeta}), \quad\mbox{where}\quad J^k(\vv{\zeta}) = \mathrm{det} \frac{\partial \vv{R}^k}{\partial \vv{\zeta}},
\end{equation}
a volume potential in $\vv{\zeta}$-reference space, $\vv{\zeta} = \left(\xi, \eta\right)^T$. In a slight abuse of notation we will when convenient use the vector notation $\vv{R}^k(\vv{\zeta})$ instead of the notation $\vv{R}^k(\xi, \eta)$ for the same function $\vv{R}^k$ (and similarly for other functions). It is known that if $\vv{R}^k$ is a $C^1$-invertible mapping onto $\widehat{\mathcal{C}}_k$, then letting $\vv{\zeta}^k_0 = \left(\vv{R}^k\right)^{-1}(\vv{r}_0)$ denote the location of a physical target $\vv{r}_0 \in \Omega$ in the reference space of the $k$\textsuperscript{th} mesh cell the identity
\begin{equation}\label{eq:Vk_hatVk_equiv}
    \mathcal{V}_k[\rho](\vv{r}_0) = \widehat{\mathcal{V}}_k[\rho \circ \vv{R}^k](\vv{\zeta}_0^k)
\end{equation}
holds~\cite[Thm.\ 5.5 \& Add.\ 5.6]{Edwards} for an integrable function $\rho$ on $\mathcal{C}_k$.

The two sections below outline (i) interpolation by a standard basis on both the straight and curved cells present in the volumetric mesh as well as (ii) the specifics of smooth quadratures and their coupling to fast algorithms. (Similar matters for the simple case of the reference box $\widehat{\mathcal{B}}_0$ are addressed in \Cref{sec:boxes}.)

\subsubsection{Smooth quadratures and interpolation on triangles}\label{sec:triangles}

This section is concerned firstly with mappings of the standard simplex into curvilinear mapped triangles and then with interpolation and quadrature of arbitrary \emph{smooth} functions over these regions; methods for singular corrections for target points $\vv{r}_0 \in \Omega$ near to or contained in $\mathcal{T}_k$ are deferred to Sections~\ref{sec:sing_quad} and~\ref{sec:sing_corr_koornwinder}.

\subsubsection*{Mapped triangles}
An arbitrary mapped triangle $\mathcal{T}_k$ can be represented using a coordinate transformation from the standard $2$-simplex
\begin{equation}\label{eq:simplex}
    \widehat{\mathcal{T}}_0 = \left\lbrace (\xi, \eta): \xi, \eta, 1 - \xi - \eta \ge 0\right\rbrace.
\end{equation}
Noting that from~\eqref{eq:Ck_Tk_Bk} the cell $\mathcal{C}_k$ and transformation $\vv{R}^k$ correspond to the cell $\mathcal{T}_k$, for $k \le N_t$, the transformation for each cell $\mathcal{T}_k$ can be written for such $k$ as
\begin{equation}\label{eq:simplex_mapping}
    \vv{R}^k(\vv{\zeta}) =  \vv{R}^k(\xi, \eta) = \begin{pmatrix}x^k(\xi, \eta)\\y^k(\xi, \eta) \end{pmatrix},\quad \vv{\zeta} = \begin{pmatrix}\xi\\ \eta\end{pmatrix},
\end{equation}
where $x^k(\xi, \eta)$ and $y^k(\xi, \eta)$ denote the $x$- and $y$-coordinate mappings from $\widehat{\mathcal{T}}_0$ to $\mathcal{T}_k$. The mappings $\vv{R}^k$ differ in the case that $\vv{R}^k$ alternatively maps into a straight-edged or a curved triangle, with the mappings $x^k$ and $y^k$ in the straight-edged triangle case given as affine maps defined such that they correctly map the coordinates of the corners of $\widehat{\mathcal{T}}_0$ to those of the three corners of $\mathcal{T}_k$. The meshing strategy described in \Cref{sec:background} also introduces curved triangles, which in our context consist specifically of curvilinear regions with two straight edges and a single curved edge. For such curved triangles $\mathcal{T}_k$ with corners $(x_1^k, y_1^k)$, $(x_2^k, y_2^k)$, and $(x_3^k, y_3^k)$ we let, without loss of generality, the curved edge connect $(x_1^k, y_1^k)$ and $(x_2^k, y_2^k)$, and introduce the transformation
\begin{equation}\label{eq:deformed_triangle_maps}
    \begin{cases}
    x^k(\xi, \eta) &= (1 - \xi - \eta) x_1^k + \xi x_2^k + \eta x_3^k + \frac{1 - \xi - \eta}{1 - \xi}\left(\lambda(\xi) - (1 - \xi) x_1^k - \xi x_2^k\right),\\
    y^k(\xi, \eta) &= (1 - \xi - \eta) y_1^k + \xi y_2^k + \eta y_3^k + \frac{1 - \xi - \eta}{1 - \xi}\left(\mu(\xi) - (1 - \xi) y_1^k - \xi y_2^k\right),\\
    \end{cases}
\end{equation}
which can easily be seen is a $C^1$-invertible map of the standard simplex $\widehat{\mathcal{T}}_0$ onto $\mathcal{T}_k$. Here, $\lambda: [0, 1] \to \mathbb{R}$ and $\mu: [0, 1] \to \mathbb{R}$ are parametrizations of the individual coordinates of the curved edge connecting $(x_1^k, y_1^k)$ and $(x_2^k, y_2^k)$ that satisfy $\lambda(0) = x_1^k$, $\lambda(1) = x_2^k$, $\mu(0) = y_1^k$, $\mu(1) = y_2^k$. This procedure for mapping $\widehat{\mathcal{T}}_0$ to an arbitrary deformed triangle $\mathcal{T}_k$ is known as the blending function method~\cite{Gordon:73a,Gordon:73b} originally introduced in the finite element literature; see also~\cite{Babuska:91}.

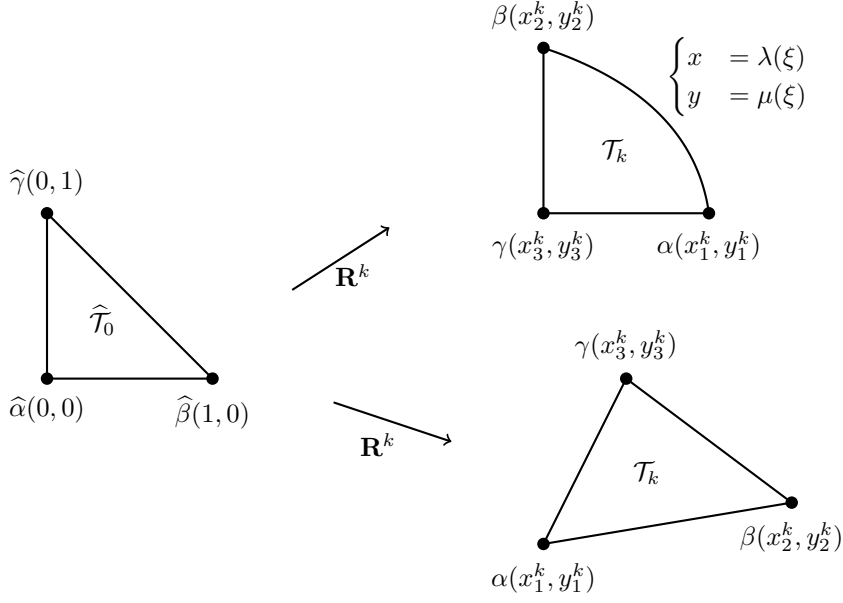
\begin{figure}
\centering
\begin{tikzpicture}[scale = 1.10]
      \draw[thick] (0, 0) -- (2, 0) {};
      \draw[thick] (0, 0) -- (0, 2) {};
      \draw[thick] (2, 0) -- (0, 2) {};
      \node[fill, circle, scale=0.5] at (0, 0) {};
      \node[anchor=north] at (0, 0-0.1) {$\widehat{\alpha}(0, 0)$};
      \node[fill, circle, scale=0.5] at (2, 0) {};
      \node[anchor=north] at (2, 0-0.1) {$\widehat{\beta}(1, 0)$};
      \node[fill, circle, scale=0.5] at (0, 2) {};
      \node[anchor=south] at (0, 2+0.1) {$\widehat{\gamma}(0, 1)$};
      \node[anchor=north] at (2/3, 3/3) {$\widehat{\mathcal{T}}_0$};
      \node (g1) at (2.85, 1) {};
      \node (g2) at (4.25, 1.9) {};
      \draw[thick, ->] (g1) to (g2);
      \node[anchor=north] at (2.7+1, 1.6-0.1) {$\mathbf{R}^k$};
      \draw[thick] (4+2, 0+2) -- (6+2, 0+2) {};
      \draw[thick] (4+2, 0+2) -- (4+2, 2+2) {};
      \draw[thick, domain=-1:1, smooth, variable=\theta, black] plot ({6+2*(1-(\theta+1)/2) + 2*2*(2*0.7-1)*(1-(\theta+1)/2)*(\theta+1)/2}, {2*(\theta+1)/2 + 2*2*(2*0.6-1)*(1-(\theta+1)/2)*(\theta+1)/2 + 2});
      \node[fill, circle, scale=0.5] at (4+2, 0+2) {};
      \node[anchor=north] at (4+2, 0+2-0.1) {$\gamma(x_3^k, y_3^k)$};
      \node[fill, circle, scale=0.5] at (6+2, 0+2) {};
      \node[anchor=north] at (6+2, 0+2-0.1) {$\alpha(x_1^k, y_1^k)$};
      \node[fill, circle, scale=0.5] at (4+2, 2+2) {};
      \node[anchor=south] at (4+2, 2+2+0.1) {$\beta(x_2^k, y_2^k)$};
      \node[anchor=north] at (2/3+4+.2+2, 3/3+2) {$\mathcal{T}_k$};
      \node[anchor=north] at (6+2.35, 1+1.25+2) {$\begin{cases} x &= \lambda(\xi) \\ y &= \mu(\xi)\end{cases}$};
      \draw[thick] (4+2, 0-2) -- (7+2, .5-2) {};
      \draw[thick] (4+2, 0-2) -- (5+2, 2-2) {};
      \draw[thick] (7+2, .5-2) -- (5+2, 2-2) {};
      \node[fill, circle, scale=0.5] at (4+2, 0-2) {};
      \node[anchor=south] at (4+2, 0-2-0.7) {$\alpha(x_1^k, y_1^k)$};
      \node[fill, circle, scale=0.5] at (7+2, .5-2) {};
      \node[anchor=south] at (7+2, .5-2-0.7) {$\beta(x_2^k, y_2^k)$};
      \node[fill, circle, scale=0.5] at (5+2, 2-2) {};
      \node[anchor=north] at (5+2, 2-2+0.7) {$\gamma(x_3^k, y_3^k)$};
      \node[anchor=north] at (4+3+0.25, 3/3-2+.1) {$\mathcal{T}_k$};
      \node (g3) at (2.35+1, -0.25) {};
      \node (g4) at (4+1, -.8) {};
      \draw[thick, ->] (g3) to (g4);
      \node[anchor=north] at (3+1, -0.55) {$\mathbf{R}^k$};
\end{tikzpicture}
    \caption{{\em\small  Depiction of reference space and mapped triangles. Left: $\widehat{\mathcal{T}}_0$ in $(\xi, \eta)$-parameter space. Right: Curved (top) and straight (bottom) mapped triangles $\mathcal{T}_k$ in physical space; note that the map $\vv{R}^k: \widehat{\mathcal{T}}_0 \to \mathcal{T}_k$ for a straight-edged triangle $\mathcal{T}_k$ takes the particularly simple form of an affine map.}}
    \label{fig:reftri_maptri}
\end{figure}

The mappings $\vv{R}^k$ are used in the method to relate integrals over triangular regions $\mathcal{T}_k$ to integrals over the standard simplex $\widehat{\mathcal{T}}_0$: since the mapping $\vv{R}^k$ given by~\eqref{eq:simplex_mapping} is a $C^1$-invertible mapping over $\widehat{\mathcal{T}}_0$ it is known~\cite[Thm.\ 5.5 \& Add.\ 5.6]{Edwards} that the integral over $\mathcal{T}_k$ can be expressed as
\begin{equation}\label{eq:Tk_cov}
    \int_{\mathcal{T}_k} \rho(\vv{r})\,\d A(\vv{r}) = \int_{\widehat{\mathcal{T}}_0} \rho(\vv{R}^k(\xi, \eta)) J^k(\xi, \eta) \,\d \xi \d \eta,
\end{equation}
for an integrable function $\rho: \mathcal{T}_k \to \mathbb{R}$, where $J^k(\xi, \eta)$ denotes the absolute value of the Jacobian determinant of the mapping $\vv{R}^k$; clearly, \Cref{eq:Tk_cov} is equivalent to~\eqref{eq:Vk_hatVk_equiv} for this class of cell. This elementary change of variables not only enables the use of an orthogonal basis for every cell $\mathcal{T}_k$, but it also enables efficient means to generate singular corrections (see \Cref{sing_precomp_quad_nodes}).

\subsubsection*{Smooth interpolation and quadrature}
At the core of the present numerical method for volume potential evaluation are robust and efficient, high-order interpolation and quadrature schemes on each of the cell types. Here we present such a scheme for the simplex using the well-known Koornwinder polynomial system. Denoting by $P_n^{(\alpha, \beta)}(\xi)$ (see~\cite[\S 22]{Abramowitz:72} for details) the Jacobi polynomial of degree $n$ that satisfies the ODE
\[
(1 - x)^2 y'' + (\beta - \alpha - (\alpha + \beta + 2)x)y' + n(n + \alpha + \beta + 1)y = 0,\quad -1 < x < 1,
\]
we will make extensive use of the Koornwinder polynomials which are defined, up to maximal total degree $p$, by
\[
    K_{nm} = \gamma_{nm} P_{n-m}^{(0, 2m+1)}(1 - 2\eta) P_m^{(0,0)}\left(\frac{2\xi}{1-\eta} - 1\right) (1 - \eta)^m; \quad m = 0, \ldots, n,\quad\mbox{and}\quad \quad n = 0, \ldots, p,
\]
where the weights $\gamma_{nm}$ are chosen so that
\[
    \int_{\widehat{\mathcal{T}}_0} K^2_{nm}(\xi, \eta) \, \d\xi\d\eta = 1.
\]
As is well-known~\cite{Koornwinder:75}, the $p(p+1)/2$ polynomials $\lbrace K_{nm}: 0 \le m \le n, \, 0 \le n < p \rbrace$ form an orthogonal basis for the space $\mathcal{P}_{p-1}$ of polynomials of total degree less than $p$ on the simplex $\widehat{\mathcal{T}}_0$. We will call an element $\rho$ of $\mathcal{P}_{p-1}$,
\[
    \rho(\xi, \eta) = \sum_{n=0}^{p-1} \sum_{m=0}^n a_{nm} K_{nm}(\xi, \eta),
\]
a $p$-th order Koornwinder expansion in keeping with classical results on the error in interpolation of smooth functions by polynomials of degree less than $p$.

The coefficients of a polynomial $\rho \in \mathcal{P}_{p-1}$ can be related to its values on a discrete set of interpolation nodes
\begin{equation}\label{eq:VR_interp_nodes}
    I_{\mathcal{T}, p} = \lbrace (\xi_{p,i}, \eta_{p,i}): 1 \le i \le N_p\rbrace, \quad N_p = p(p+1)/2.
\end{equation}
In detail, denoting by $\vv{a}$ and $\vv{\rho}$ the vectors with elements $(\vv{a})_{nm} = a_{nm}$ and $(\vv{\rho})_{i} = \rho(\xi_{p,i}, \eta_{p,i})$ the coefficients will satisfy
\begin{equation}\label{eq:vals_to_coeffs}
    \vv{V}_p \vv{a} = \vv{\rho},
\end{equation}
where $(\vv{V}_p)_{nm, i} = K_{nm}(\xi_{p,i}, \eta_{p,i})$ is the so-called coefficients-to-values map associated with the Koornwinder polynomials on the nodal set $I_{\mathcal{T}, p}$. Provided the nodal set $I_{\mathcal{T}, p}$ is such that $\vv{V}_p$ is nonsingular we denote the values-to-coefficients map $\vv{C}_p = \vv{V}_p^{-1}$---while in our case $\vv{V}_p$ will always be invertible due to a specific choice of nodal set $I_{\mathcal{T}, p}$, see also~\cite{Chung:77} for a sufficient condition for a generic nodal set $I_{\mathcal{T}, p}$ to yield an invertible matrix $\vv{V}_p$.

Recent contributions~\cite{Vioreanu:14} have developed nodal sets $I_{\mathcal{T}, p}$ leading to favorable conditioning of the maps $\vv{V}_p$ and $\vv{C}_p$ that we utilize, and, making this selection we henceforth denote by $I_{\mathcal{T}, p}$ the Vioreanu-Rokhlin nodes for interpolation by polynomials with total degree less than $p$. It will also be convenient to introduce the oversampling matrix $\vv{O}_{p_1,p_2}$ which maps values at the Vioreanu-Rokhlin nodes $I_{\mathcal{T}, p_1}$ in a Koornwinder expansion of order $p_1$ to that same expansions' values at the Vioreanu-Rokhlin nodes $I_{\mathcal{T}, p_2}$ for an expansion of order $p_2 \ge p_1$.
For a generic smooth function $\rho: \widehat{\mathcal{T}}_0 \to \mathbb{R}$ we have  the approximation
\begin{equation}
    \rho(\xi, \eta) \approx \rho_p(\xi, \eta) = \sum_{n=0}^{p-1} \sum_{m=0}^n a_{nm} K_{nm}(\xi, \eta),
\end{equation}
accurate to $p$-th order. The expansion coefficients $a_{nm}$ in a $p$-th order Koornwinder expansion are determined for generic smooth functions and polynomials alike via the solution to the system~\eqref{eq:vals_to_coeffs}.

Having considered interpolation, we turn to quadratures of smooth functions: a target $\vv{r}_0$-independent high-order quadrature rule is required on the simplex for integrals such as those in \Cref{eq:Tk_cov} with $\rho$ smooth. (The $\vv{r}_0$-independence of the rule---that the source nodes are independent of the target---is required for compatibility with fast summation techniques.) For this task, we turn to generalized Gaussian quadrature rules, in which context we recall that one-dimensional Gaussian quadrature yields a specific set $\{(\xi_i, w_i);\, i = 1, 2, \ldots, N\}$ of $N$ quadrature nodes and weights results that can integrate exactly all polynomials of degree at most $2N-1$ (the weights $w_i$ following from a specific choice of nodal set). In higher dimensions, it appears that perfect Gaussian quadratures are unfortunately not available, but \emph{generalized} Gaussian quadratures have been introduced~\cite{XiaoGimbutas:10,Vioreanu:14} which, for $N$ quadrature nodes in $d > 1$ dimensions, rather than exactly integrating $dN$ functions (as would be the case for a perfect Gaussian quadrature rule), instead only integrate some number of functions greater than $N$. Fortunately, the interpolation nodes $I_{\mathcal{T}, p}$ introduced above also have associated with them corresponding weights, thereby furnishing us with a highly-efficient set of quadrature nodes and weights (the \emph{efficiency} of a Gaussian quadrature is the ratio of the number of functions integrated exactly to the ideal ``Gaussian'' number $dN$)---for details see~\cite{Vioreanu:14}. Thus, for a smooth function $\rho$ on the domain $\mathcal{T}_k$ and in view of~\eqref{eq:Tk_cov} we use the interpolatory quadrature rule
\begin{equation}\label{eq:vio_rokh_quad_rule}
    \int_{\mathcal{T}_k} \rho(\vv{r}) dA(\vv{r}) \approx \sum_{j=1}^{N_p} \rho\left(\vv{R}^k(\xi_j, \eta_j)\right) J^k(\xi_j, \eta_j) w_j, \quad\mbox{where}\quad (\xi_j, \eta_j) \in I_{\mathcal{T}, p},
\end{equation}
termed so because the quadrature nodes $(\xi_i, \eta_i)$ coincide with the set $I_{\mathcal{T},p}$ of interpolation nodes; it should be cautioned that this is not a quadrature rule of order $p$~\cite[\S 5]{Vioreanu:14}.

\subsubsection*{Oversampled smooth quadratures}

An ideal smooth quadrature rule for $\mathcal{V}_k[f](\vv{r}_0)$ delivers accurate approximations with error on the level of interpolation of the function $f$ over $\mathcal{T}_k$, but, unfortunately, large gradients in the integrand due to $G$ over cells $\mathcal{T}_k$ near to $\vv{r}_0$ can lead to quadrature error that dominates that of interpolation of $f$ at any fixed order $p$. Oversampled quadratures can address this problem, limiting the number of nodal points at which $f$ is required while delivering high accuracy, and we describe next their use. Defining the vector of function samples $\vv{f}^k_{p}$ by
\[
    (\vv{f}_{p}^k)_i = f\left(\vv{R}^k(\xi_{p,i}, \eta_{p,i})\right), \quad\mbox{where}\quad (\xi_{p,i}, \eta_{p, i}) \in I_{\mathcal{T}, p}\mbox{ for } i = 1, \ldots, N_p,
\]
so that
\[
    f(\vv{R}^k(\xi, \eta)) \approx \sum_{n=0}^{p-1} \sum_{m=0}^n a^k_{nm} K_{nm}(\xi, \eta) \quad\mbox{with}\quad \vv{V}_p \vv{a}^k = \vv{f}_p^k,
\]
interpolated approximate values of the source $f$ are obtained at nodes $I_{\mathcal{T}, q}$ via the oversampling map $\vv{O}_{p_1,p_2}$ with $p_1 = p$ and $p_2 = q \ge p$ (see the previous section on smooth interpolation),
\[
    \vv{f}_{q}^k = \vv{O}_{p,q} \vv{f}_p^k.
\]
where the integer $q$ represents the degree of oversampled quadrature. The quadrature rule~\eqref{eq:vio_rokh_quad_rule} applied to $\rho(\vv{r}) = G(\vv{r}, \vv{r}_0) f(\vv{r})$ for each $\mathcal{T}_k$ satisfying $\mathcal{T}_k \subset \mathcal{C}^\mathrm{far}(\vv{r}_0)$ is then
\begin{equation}\label{eq:Tk_smooth_quad}
    \mathcal{V}_k[f](\vv{r}_0) \approx \sum_{j=1}^{N_{q}} G(\vv{R}^k(\xi_j, \eta_j), \vv{r}_0) (\vv{f}_{q}^k)_j J^k(\xi_j, \eta_j) w_j.
\end{equation}
One detail to note is that while quadrature occurs in reference $(\xi, \eta)$ space in which context the function is not immediately obviously amenable to standard fast algorithms, the sum in~\eqref{eq:Tk_smooth_quad} is nevertheless compatible with FMMs in physical space by viewing the quantity $J^k(\xi_j, \eta_j) w_j$ as modified weights for the discrete inner product.

\subsection{A singular quadrature scheme for mapped simplices}\label{sec:sing_quad}
This section is devoted to a description of a singular quadrature scheme that can evaluate volumetric integrals over arbitrary star-shaped regions of the plane. While the singular quadrature methods outlined below are quite generic geometrically (though see \Cref{rem:starshaped_vs_convex}), for our purposes here singular integral evaluation will be required over one of precisely two volumetric regions $\widehat{\mathcal{C}}$: a box $\widehat{\mathcal{C}} = \widehat{\mathcal{C}}_k = \widehat{\mathcal{B}}_0$ and a standard simplex $\widehat{\mathcal{C}} = \widehat{\mathcal{C}}_k = \widehat{\mathcal{T}}_0$ (as is suggested by the notation, we suppress the dependence on $k$ in this section whenever possible). This section begins by presenting the basic theory of converting convolutions of the Green's function with a source function $\rho$ over certain star-shaped domains to integrals over each of those domain's boundaries (using ideas related to Poincar\'e's lemma), then discusses treatment of the singular kernel in \Cref{sec:singquad} and finally addresses some subtler geometric details in \Cref{sec:close}. In principle the function $\rho$ could incorporate directly the actual source function $f$ (with mapped argument) arising from the PDE boundary value problem \Cref{eq:ellipt_prob_pde} but in practice will be chosen to be related to a member of a given family of orthogonal polynomials, as detailed in \Cref{sec:sparse_func_approx}.

\subsubsection{Poincar\'e's lemma and volume-to-boundary integral conversion for $\widehat{\mathcal{V}}_k$}
It is first useful to recall some terminology. A region $\mathcal{S}$ is called star-shaped if there exists \emph{some} point $\vv{r}_* \in \mathcal{S}$ such that every line segment that connects $\vv{r}_*$ to any other point $\vv{r} \in \mathcal{S}$ \emph{lays entirely in $\mathcal{S}$} (the region is called star-shaped with respect to $\vv{r}_*$). On the other hand, a star-shaped region is convex if and only if it is star-shaped with respect to every $\vv{r}_* \in \mathcal{S}$.

Inspired by elements of the proof of Poincar\'e's lemma, which utilizes certain maps involving iterated boundary integrals, we make use of an extension of these ideas for integrals involving singular functions, but first state its main holding in a simpler lemma limited to smooth functions. A proof is given in the appendix, using the perhaps more familiar tools of vector calculus (cf.\ elements of this lemma using the language of differential forms in~\cite{spivak2018calculus}).

\begin{lemma}\label{lem:smooth_poincare}
    Let the closed region $\mathcal{S} \subset \mathbb{R}^2$ be star-shaped with respect to the origin with a piecewise smooth boundary $\partial \mathcal{S}$ and assume $\rho: \mathcal{S} \to \mathbb{R}$ is continuously differentiable on $\mathcal{S}$. Then the relation
    \begin{equation}\label{eq:smooth_poincare}
        \int_{\mathcal{S}} \rho(\vv{r})\,\mathrm{dA} = \oint_{\partial \mathcal{S}}\left(\int_0^1 t\,\rho(t\vv{r})\,\mathrm{d}t\right)  \vv{r}\times \vv{\tau}\,\mathrm{d}s
    \end{equation}
    holds, where $s$ denotes the arclength, $\vv{\tau}$ denotes the unit tangent vector of $\partial \mathcal{S}$, and $\vv{r}\times \vv{\tau}$ is understood as a scalar.
\end{lemma}

The integrals resulting from \Cref{lem:smooth_poincare} have at times been referred to as dilation integrals (with dilation with respect to the star point), with e.g.\ reference~\cite{Barnett:21} relying precisely (up to a fixed geometry-dependent translation) on the relation~\eqref{eq:smooth_poincare} to generate accurate quadratures of smooth functions. However, in the context of volume potentials with weakly singular kernels, such approaches~\cite{Fata:12,Gao:02,Gao:05,Wang:17a,Wang:17b} have previously faced a number of challenges related to star-shapedness, location of quadrature points, and presence of singularities (see the discussion in the Introduction); in some cases~\cite{Fata:12} the applicability is restricted to Poisson problems.

\begin{remark}\label{rem:starshaped_vs_convex}
    For a \emph{convex} reference cell $\mathcal{S} = \widehat{\mathcal{C}}$, which in the proposed method are the only domains over which quadrature is performed, \Cref{lem:poincare} below applies to \emph{every} point $\vv{\zeta}_* \in \widehat{\mathcal{C}}$, while this may unfortunately not hold for a generic star-shaped region. It is a significant strength of the present mapped-domain approach, in contrast to the possibility of applying this lemma in physical space (for which case difficulties may arise when targets lay in certain subsets of nonconvex cells), that the Poincar\'e lemma-related ideas can be successfully applied for every point $\vv{\zeta}_*$. More generally, methods based on transformation of volume integrals to surface integrals potentially suffer when integration domains are non-convex as they generally result in quadrature points laying outside the integration domain result~\cite{Gunderman:21,Fata:12,Muller:13,Gao:02,Gao:05}, a long-recognized issue in multidimensional quadrature~\cite{Cools:97,Maxwell:87}. As a more minor matter, the presence of such quadrature nodes implies the existence of negative quadrature weights, which is generally considered unfavorable in view of stability concerns.
\end{remark}

The following corollary of \Cref{lem:smooth_poincare} is used in what follows to express the volume potentials in terms of boundary integrals; its proof is given in the appendix.

\begin{corollary}\label{lem:poincare}
    Let $K = K(\vv{r})$ denote a weakly-singular kernel function which is continuously differentiable for $\vv{r}\in \mathbb{R}^2\setminus \{\vv{0}\}$, assume the closed region $\mathcal{S} \subset \mathbb{R}^2$ has a piecewise smooth boundary $\partial\mathcal{S}$, and assume the function $\rho: \mathcal{S} \to \mathbb{R}$ is continuously differentiable on $\mathcal{S}$. Then for each $\vv{r}_* \in \mathcal{S}$ such that $\mathcal{S}$ is star-shaped with respect to $\vv{r}_*$, we have
    \begin{equation}\label{eq:lemma1}
        \int_{\mathcal{S}} K(\vv{r}-\vv{r}_*) \rho(\vv{r}) \,\mathrm{dA} = \oint_{\partial \mathcal{S}}\left(\int_0^1 t\,K(t\left(\vv{r}-\vv{r}_*\right)) \rho(t(\vv{r} - \vv{r}_*) + \vv{r}_*)\,\mathrm{d}t\right)  \left((\vv{r} - \vv{r}_*)\times \vv{\tau}\right)\,\mathrm{d}s.
    \end{equation}
\end{corollary}

\begin{figure}[!ht]
    \centering
    \includegraphics[height=0.33\textwidth]{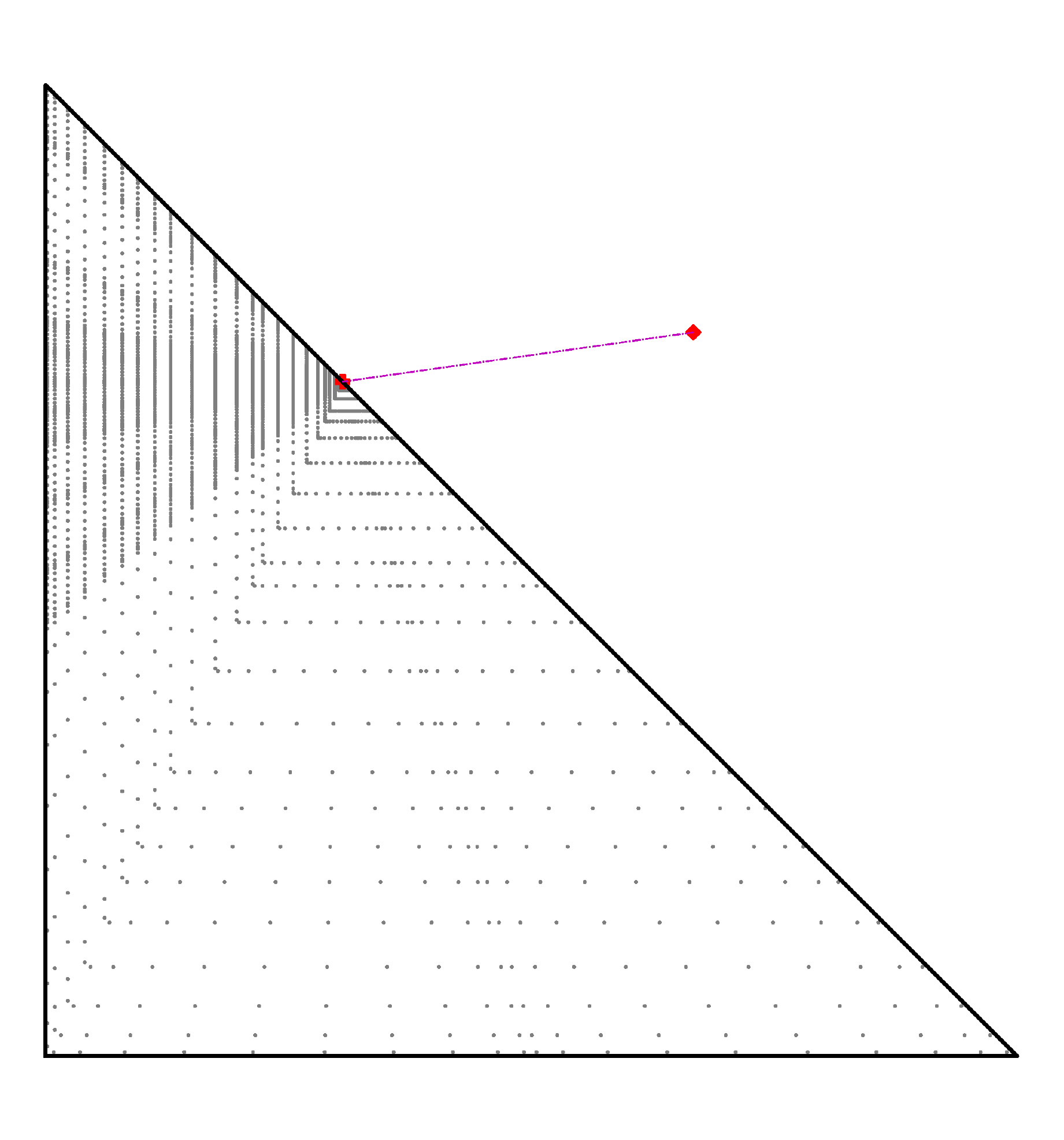}
    \includegraphics[height=0.33\textwidth]{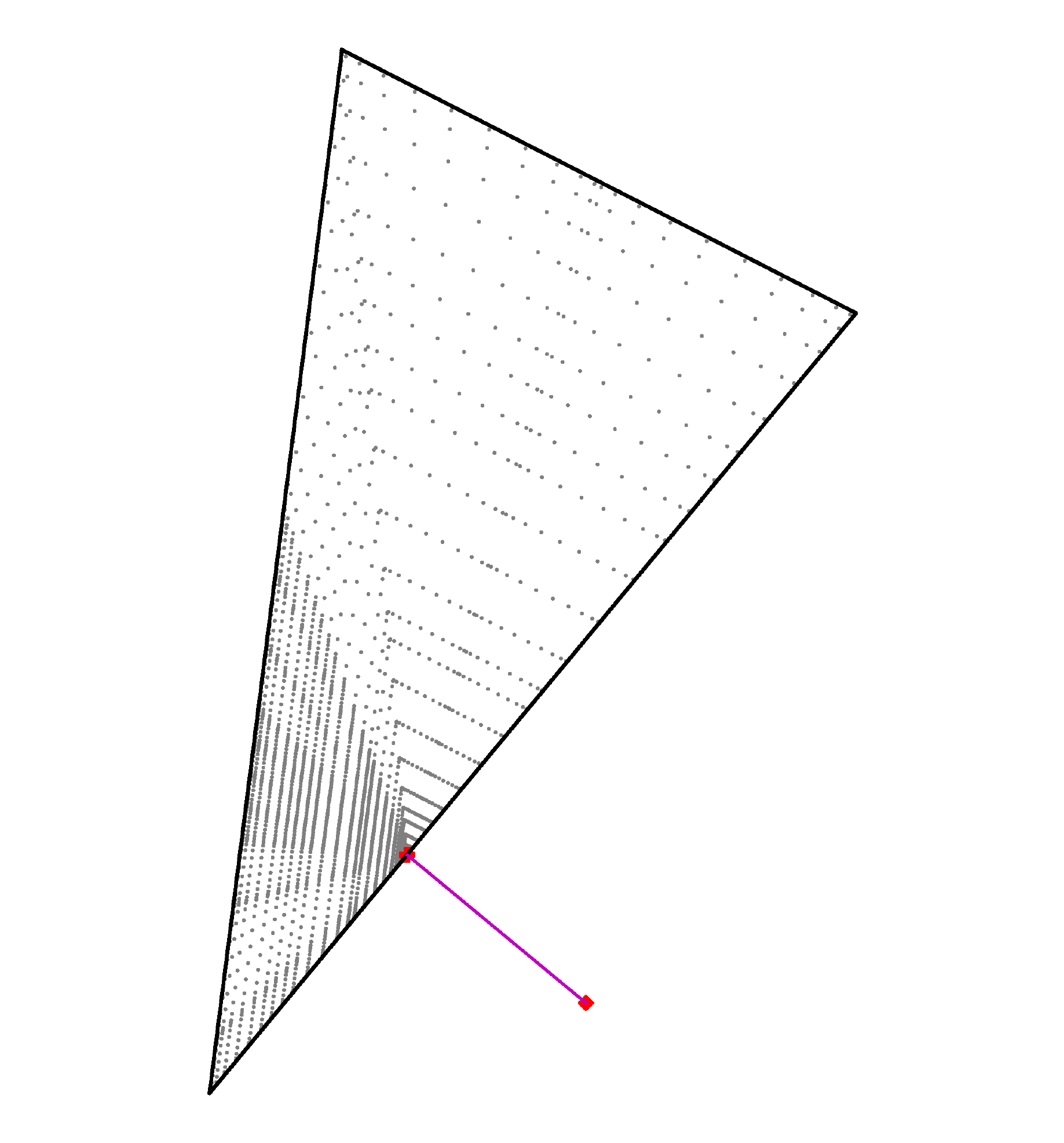}
    \caption{{\em\small In reference (left) and physical (right) space for a straight triangle cell, the high-order quadrature nodes (arising from a $10$\textsuperscript{th}-order modified Gaussian rule for the $I^k$ integral and a $10$\textsuperscript{th}-order composite Gauss--Legendre rule for the integral over $\partial\widehat{\mathcal{C}}$) are displayed that result from the proposed quadrature scheme applied to a near-singular target point $\vv{r}_0 \not\in \mathcal{C}$ (marked in red). The dashed line segment connects the target $\vv{r}_0$ to the nearest physical point $\vv{r}_* \in \mathcal{C}$ which is selected as the star-point (also marked in red); in reference domain $\vv{\zeta}_* = \left(\vv{R}\right)^{-1}(\vv{r}_*)$ is \emph{not} the closest point to $\vv{\zeta}_0$ since the map $\vv{R}$ is affine (see also \Cref{eq:star_point}).}}
    \label{fig:nearsingquad_nodes}
\end{figure}

The proposed methodology identifies for every singular and near-singular target $\vv{r}_0 \in \Omega$ a point $\vv{r}_* = \vv{r}_*(\vv{r}_0) \in \mathcal{C}_k$ using the rule
\begin{equation}\label{eq:rstar_def}
    \vv{\zeta}_* = \argmin_{\vv{\zeta} \in \widehat{\mathcal{C}}_k} \left|\vv{R}^k(\vv{\zeta}) - \vv{r}_0\right|\quad\mbox{and}\quad \vv{r}_* = \vv{R}^k(\vv{\zeta}_*).
\end{equation}
Having identified $\vv{\zeta}_* \in \widehat{\mathcal{C}}_k$ we now apply \Cref{lem:poincare} and then describe quadrature rules for evaluating certain resulting integrals in the reference domain $\widehat{\mathcal{C}} = \widehat{\mathcal{C}}_k$.
\begin{remark}\label{eq:star_point}
  A subtle but important point for near-singular targets $\vv{r}_0$ is that the star point $\vv{r}_*$ must be selected as the closest point in $\partial\mathcal{C}_k$ to $\vv{r}_0$ despite quadrature occurring over $\widehat{\mathcal{C}}$, as this point corresponds to the minimum distance $d$ that arises as an argument to the singular kernel (cf.~\Cref{eq:near_sing_dist_rule} and \Cref{fig:nearsingquad_nodes}). By adapting to this point the quadrature rule is able to deliver optimal accuracy.
\end{remark}

Writing $\widehat{\mathcal{V}}_k[\rho]$ in the form
\[
    \widehat{\mathcal{V}}_k[\rho] = \int_{\widehat{\mathcal{C}}} G\left(\vv{R}^k(\vv{\zeta} - \vv{\zeta}_* + \vv{\zeta}_*) - \vv{R}^k(\vv{\zeta}_0)\right) J^k(\vv{\zeta}) \rho(\zeta)\,\mathrm{d}A
\]
and making the selections of region $\mathcal{S} = \widehat{\mathcal{C}}$, kernel $K(\vv{\tau}) = G\left(\vv{R}^k(\vv{\tau} + \vv{\zeta}_*) - \vv{R}^k(\vv{\zeta}_0)\right)$ and smooth source $J^k(\vv{\zeta}) \rho(\vv{\zeta})$ in \Cref{lem:poincare}, we obtain
\begin{equation}\label{eq:poincare}
\begin{split}
    \widehat{\mathcal{V}}_k[\rho]\left(\vv{\zeta}_0\right) = \oint_{\partial \widehat{\mathcal{C}}} \left( \int_0^1 t \right.&\left.\vphantom{\int_0^1}G(\vv{R}^k(\vv{\zeta}_*+t\left(\vv{\zeta}-\vv{\zeta}_*\right)),\vv{R}^k(\vv{\zeta}_0))\right.\times\\
    & \times \left.\vphantom{\int_0^1}J^k( \vv{\zeta}_*+t\left(\vv{\zeta}-\vv{\zeta}_*\right))\rho(\vv{\zeta}_*+t\left(\vv{\zeta}-\vv{\zeta}_*\right))\,\mathrm{d}t\right) \left((\vv{\zeta} - \vv{\zeta}_*)\times\vv{\tau}\right)\,\mathrm{d}s.
\end{split}
\end{equation}
The expression~\eqref{eq:poincare} for the volume potential $\widehat{\mathcal{V}}_k[\rho]$ can be written in the simplified form
\begin{equation}\label{eq:poincare_Ik}
    \widehat{\mathcal{V}}_k[\rho]\left(\vv{\zeta}_0\right) = \oint_{\partial \widehat{\mathcal{C}}} I^k[\rho](\vv{\zeta}, \vv{\zeta}_0) \left((\vv{\zeta} - \vv{\zeta}_*)\times \vv{\tau}\right)\,\mathrm{d}s,
\end{equation}
where $I^k[\rho](\vv{\zeta}, \vv{\zeta}_0)$ is defined as the inner integral of the iterated integral:
\begin{equation}\label{eq:Ik_def}
    I^k[\rho](\vv{\zeta}, \vv{\zeta}_0) = \int_0^1 t\,G(\vv{R}^k(\vv{\zeta}_*+t\left(\vv{\zeta}-\vv{\zeta}_*\right)),\vv{R}^k(\vv{\zeta}_0))J^k( \vv{\zeta}_*+t\left(\vv{\zeta}-\vv{\zeta}_*\right)) \rho(\vv{\zeta}_*+t\left(\vv{\zeta}-\vv{\zeta}_*\right))\,\mathrm{d}t.
\end{equation}
The presence of the Green function (near-)singularity at $t=0$, for all $\vv{\zeta}, \vv{\zeta}_0$, in the one-dimensional integral in~\eqref{eq:Ik_def} suggests that one-dimensional quadrature schemes for the integral $I^k[\rho](\vv{\zeta}, \vv{\zeta}_0)$ may be effective, which we outline in the next section. Before discussing these quadrature schemes we first introduce the quadrature scheme for the outer integral in~\eqref{eq:poincare_Ik}.

Since $\partial\widehat{\mathcal{C}}$ is piecewise smooth, it follows that $I^k[\rho](\vv{\zeta}, \vv{\zeta_0})$ is piecewise smooth as a function of $\vv{\zeta} \in \partial \widehat{\mathcal{C}}$, implying that the integrand of the circulation integral \Cref{eq:poincare_Ik} is likewise piecewise smooth. The quadrature rule we propose for $\widehat{\mathcal{V}}_k[\rho](\vv{\zeta}_0)$ is developed by segmenting the boundary $\partial\widehat{\mathcal{C}}$ into $M_{\partial\widehat{\mathcal{C}}}$ intervals, and obtaining a composite Gauss-Legendre rule by application of a $P$-point Gauss--Legendre quadrature rule on each interval, with $P$ a fixed integer (we use the selection $P = 16$ for the numerical examples in this article). We thus have the quadrature rule
\begin{equation}\label{eq:Vk_bdry_quad}
    \widehat{\mathcal{V}}_k[\rho]\left(\vv{\zeta}_0\right) \approx \sum_{j=1}^{PM_{\partial\widehat{\mathcal{C}}}} I^k[\rho](\vv{\zeta}_j, \vv{\zeta}_0) \left((\vv{\zeta}_j - \vv{\zeta}_*)\times\vv{\tau}_j\right) w_j,
\end{equation}
where $\{\vv{\zeta}_j\}_{j=1}^{PM_{\partial\widehat{\mathcal{C}}}}\subset \partial\widehat{\mathcal{C}}$ and $\{w_j\}_{j=1}^{PM_{\partial\widehat{\mathcal{C}}}}$ are the quadrature nodes and weights. The error arising from use of this quadrature can depend in an important manner on the proximity of $\vv{\zeta}_0$ relative to the boundary $\partial \widehat{\mathcal{C}}$, and, in particular, we find that a uniform distribution of composite Gauss-Legendre intervals is not sufficient to achieve desired accuracies. Further details on boundary discretization are given in \Cref{sec:close} where an optimal interval distribution is described which may depend on the target $\vv{\zeta}_0$.

\subsubsection{Singular and near-singular quadrature rules for $I^k[\rho]$}\label{sec:singquad}
While use of \Cref{lem:poincare} transforms the volume integral over the cell $\widehat{\mathcal{C}}$ to an iterated integral over the boundary of the cell $\partial \widehat{\mathcal{C}}$, the integrand of the inner integral remains potentially non-smooth (indeed, potentially singular depending on the strength of the kernel singularity). The one-dimensional quadrature scheme outlined in what follows is tailored to the known singular behavior of the kernel: crucially, only the endpoint \emph{asymptotic behavior} of the integrand is relevant to quadrature rule selection, making the proposed scheme applicable to many kernel functions, including, but not limited to, those arising in a variety of Green functions of mathematical physics (in particular, the kernel function could be more singular than the Green functions typically encountered in elliptic PDEs).
\begin{remark}\label{rem:singular}
    For convenience, we refer to the one-dimensional integrals $I^k[\rho]$ as being singular (for $\vv{\zeta}_0 \in \widehat{\mathcal{C}}_k$) or near-singular (for $\vv{\zeta}_0 \not\in \widehat{\mathcal{C}}_k$), even though, in view of the $t$-factor present in the integrand of $I^k[\rho]$, it is possible that the integrand is not truly singular; for example, for the elliptic PDEs we consider in this article, the small-$t$ asymptotic behavior of $G$ in $I^k$ is merely logarithmic ($s(t) = \log t$ in~\eqref{eq:alpert_func} and $S(t) = \log t$ in~\eqref{eq:kolmrokh_int} below) and the integrand of $I^k[\rho]$ is thus merely non-smooth. Indeed the $t$-factor, which could be loosely viewed as an analogue of weights arising in polar changes of variables, allows the treatment of highly singular kernels by the proposed methodology. In any case, for the PDE kernels considered in the examples of this article, singular and near-singular quadrature rules are still required for high-order accuracy (i.e.\ the $t$-factor is implicitly treated as part of the smooth component of the integrand ($\phi(t)$ in~\eqref{eq:alpert_func} and $k_2(t)$ in~\eqref{eq:kolmrokh_int} below)).
\end{remark}

\begin{figure}[!ht]
    \centering
    \includegraphics[height=0.33\textwidth]{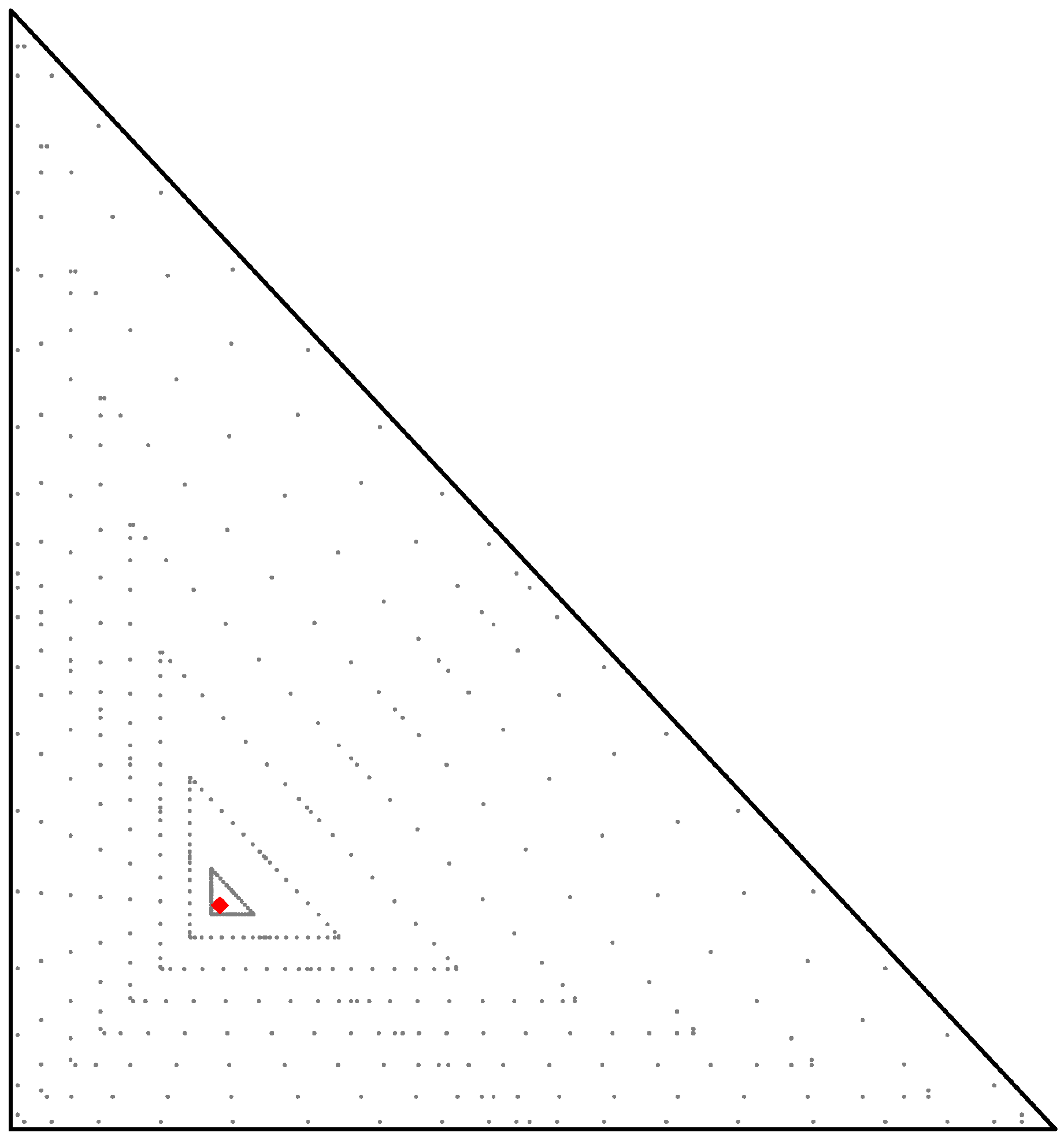}
    \includegraphics[height=0.33\textwidth]{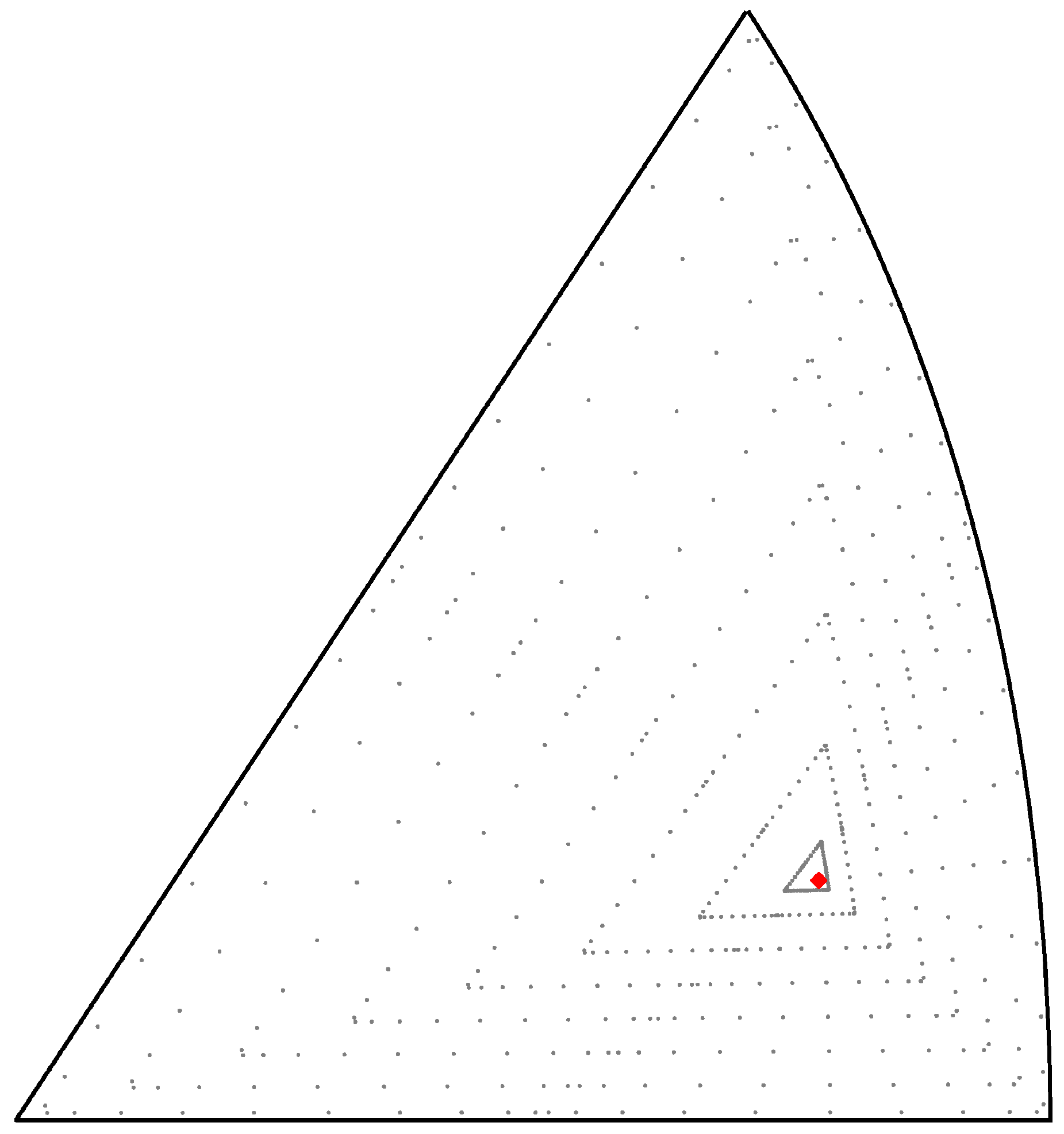}
    \caption{{\em\small In reference (left) and physical (right) space for a curvilinear cell, the high-order quadrature nodes (arising from a $5$\textsuperscript{th}-order corrected trapezoidal rule for the $I^k$ integral and a $10$\textsuperscript{th}-order composite Gauss--Legendre rule for the integral over $\partial\widehat{\mathcal{C}}$) are displayed that result from the proposed quadrature scheme applied to a singular target point $\vv{r}_0 \in \mathcal{C}$ (marked in red).}}
    \label{fig:singquad_nodes}
\end{figure}

In detail, the primary difficulty in accurate evaluation of the $I^k[\rho](\vv{\zeta}, \vv{\zeta}_0)$ integral is that, in view of the nature of the integrand in $I^k$, the target point $\vv{\zeta}_0$ will lie precisely at (in the case $\vv{\zeta}_0 \in \widehat{\mathcal{C}}$) or instead near (in the case $\vv{\zeta}_0 \not\in \widehat{\mathcal{C}}$) one end of the line segment connecting $\vv{\zeta}_*$ ($t = 0$) and $\vv{\zeta} \in \partial\widehat{\mathcal{C}}$ ($t = 1$). In the singular case, there is possibly an integrable singularity in the integrand at the left ($t=0$) endpoint of the integration interval due to the singular nature of the Green function $G$. Problems remain even in the case when there is a positive distance of $\vv{\zeta}_0$ to $\widehat{\mathcal{C}}$, since steep gradients, again at $t=0$, can be challenging to accurately resolve with any fixed target-independent quadrature scheme. The high-order accurate quadrature scheme described in this article depends on a corresponding high-order accurate quadrature rule for the integral $I^k$, of which a wide variety of suitable schemes for singularity behavior of various types have been developed over many years, e.g.~\cite{Alpert:99,Kapur:97,Kolm:01} (see~\cite{Rokhlin:90} and references therein for discussion of early work in this direction). The integrals $I^k$, when the integrand is non-smooth, are amenable to the use of existing corrected trapezoidal quadrature rules for functions with known (singular) endpoint behavior and we utilize in this article the rules of reference~\cite{Alpert:99} for singular quadrature---rules for integrating smooth functions multiplied by singular functions of logarithmic and (integrable) inverse-power type; for near-singular quadrature we use the rules introduced in reference~\cite{Kolm:01}.

We first detail the application of the Alpert~\cite{Alpert:99} quadrature rule to $I^k$ for self-interaction (singular) terms. The Alpert rules provide endpoint-corrected trapezoidal quadrature nodes and weights for integrals of functions $h(t): (0, 1] \to \mathbb{R}$ of the form
\begin{equation}\label{eq:alpert_func}
    h(t) = \phi(t)s(t) + \psi(t),
\end{equation}
 where $\phi(t), \psi(t) \in C^k[0, 1]$ and $s(t) \in C(0, 1]$ is an integrable function that is singular at $t=0$. Letting the set of pairs $\{(t_i,v_i),\, 1\leq i\leq m_P \}$ denote, for a given integer $P$ (the selection $P = 16$ is made in all numerical results in this article), the nodes and weights of a $P$-th order convergent Alpert quadrature rule for a given $s(t)$, we obtain from~\eqref{eq:Ik_def}
\begin{equation}\label{eq:Ik_vp_singular}
    I^k[\rho]\left(\vv{\zeta}, \vv{\zeta}_0\right) \approx \sum_{i = 1}^{m_P} t_i\, G\left(\vv{R}^k(\vv{\zeta}_*+t_i\left(\vv{\zeta}-\vv{\zeta}_*\right)),\vv{R}^k(\vv{\zeta}_0)\right)J^k( \vv{\zeta}_*+t\left(\vv{\zeta}-\vv{\zeta}_*\right)) \rho\left(\vv{\zeta}_*+t_i \left(\vv{\zeta}-\vv{\zeta}_*\right)\right)v_i,
\end{equation}
which is a $P$-th order approximation to $I^k[g]$ and where we used the selection  $\vv{\zeta}_* = \vv{\zeta}_0$ that has been made for the case of singular target points ($\vv{\zeta}_0 \in \widehat{\mathcal{C}}$).

In the near-singular case ($\vv{\zeta}_0 \not\in \widehat{\mathcal{C}}$) we turn to modified Gaussian quadrature rules for integrands with known singular behavior $S(t)$ as $t \to 0^+$. Such methods provide, for a given $d$ which lays in intervals of the form $[10^{-q-1}, 10^{-q}]$ for positive integer $q$, quadrature nodes and weights $\lbrace (\tilde{t}_i, \tilde{v}_i), 1 \le i \le m'_P \rbrace$ so that the near-singular rule
\begin{equation}\label{eq:kolmrokh_int}
\int_0^1 \left( k_1(t) + k_2(t) S(t + d) \right)\,\d t \approx \sum_{i=1}^{m'_p} \left(k_1(\tilde{t}_i) + k_2(\tilde{t}_i) S(\tilde{t}_i + d)\right) \tilde{v}_i,
\end{equation}
holds to high accuracy (to within an accuracy $\varepsilon$ of $\varepsilon \approx 10^{-15}$) for $k_1$ and $k_2$ polynomials of degree at most $P$ (were the quadrature rule exact for such polynomials it would be a perfect Gaussian, hence the `modified' moniker). Methods and theory are described in~\cite{Kolm:01} for generation of such quadratures for a variety of singular behaviors; for the two-dimensional elliptic PDE demonstrations of the present work, we rely on pre-computed rules~\cite{Hao:14} for the case $S(t) = \log(t)$, where for the selection $P = 10$ made everywhere in this article a total of $m'_P = 24$ nodes are required. These quadrature rules are parametrized by the ``near-singularity distance'' $d$ defined by the expression
\begin{equation}\label{eq:near_sing_dist_rule}
    d = \left| \vv{R}^k(\vv{\zeta}_*) - \vv{R}^k(\vv{\zeta}_0)\right|,
\end{equation}
which is the minimum distance that can occur for arguments to the kernel for this integration domain and target $\vv{r}_0$ (cf.~\eqref{eq:rstar_def}). We thus have the quadrature rule for $I^k[\rho](\vv{\zeta}, \vv{\zeta}_0)$,
\begin{equation}\label{eq:Ik_vp_nearsingular}
    I^k[\rho](\vv{\zeta}, \vv{\zeta}_0) \approx \sum_{i = 1}^{m'_P} \tilde{t}_i\, G\left(\vv{R}^k\left(\vv{\zeta}_*+\tilde{t}_i\left(\vv{\zeta}-\vv{\zeta}_*\right)\right),\vv{R}^k\left(\vv{\zeta}_{0}\right)\right) J^k( \vv{\zeta}_*+t\left(\vv{\zeta}-\vv{\zeta}_*\right)) \rho\left(\vv{\zeta}_*+\tilde{t}_i \left(\vv{\zeta}-\vv{\zeta}_*\right)\right)\tilde{v}_i.
\end{equation}

\subsubsection{Close evaluation}\label{sec:close}
\begin{figure}[!ht]
    \centering
    \includegraphics[width=0.305\textwidth]{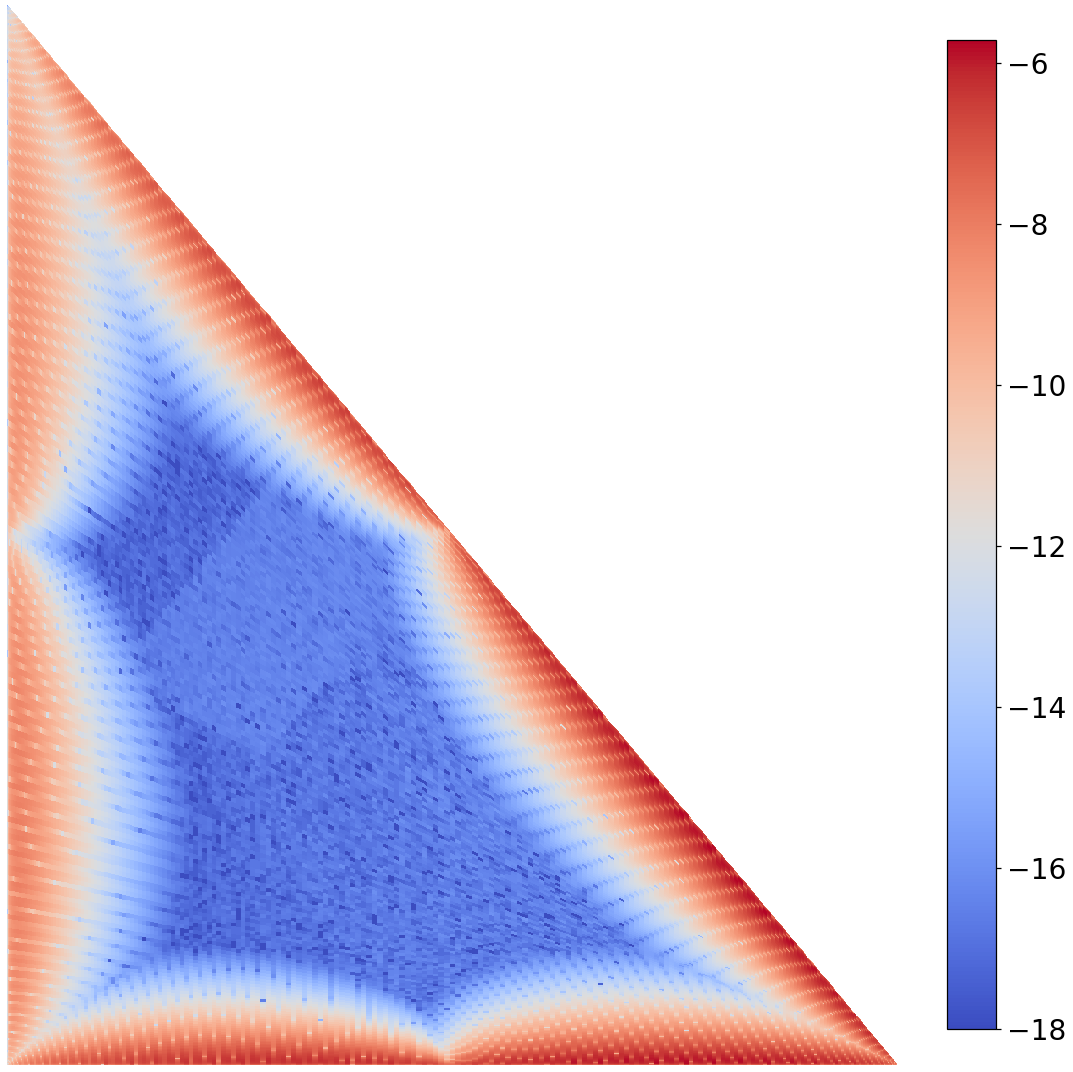}
    \quad\quad\quad\quad
    \includegraphics[width=0.33\textwidth]{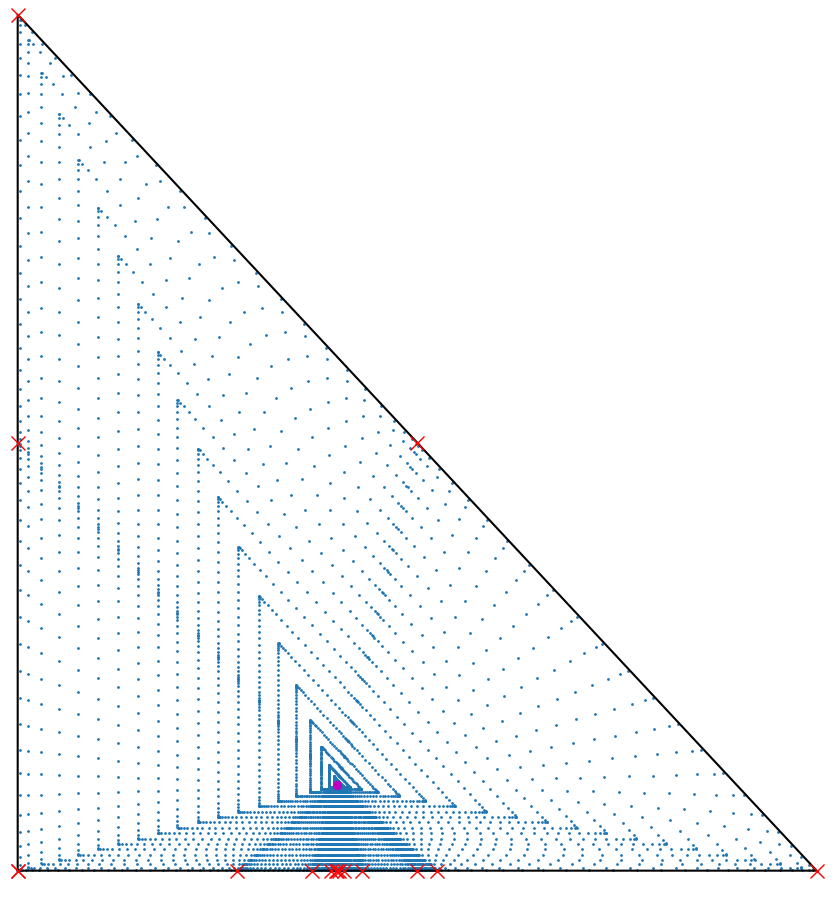}
    \caption{{\em\small Left: Additional error (in evaluation of a logarithmic-kernel volume potential) that arises from a (na\"ive) uniform distribution of intervals for purposes of a boundary discretization in \Cref{eq:poincare_Ik}, for each target point in the interior of $\widehat{\mathcal{C}}$ a unit simplex (demonstrating the need for some special care in handling target points laying close to a boundary); the color indicates the base-$10$ logarithm of the `error', compared to the results of the proposed non-uniform boundary interval distribution that is \emph{target-adapted}. Right: Singular volumetric quadrature nodes (blue) and boundary interval endpoints (red) resulting from the proposed target-adapted methodology when the target (magenta) lays close to the boundary of the reference cell $\widehat{\mathcal{C}}$; the proposed boundary interval distribution produces uniformly-accurate volume potential values over the unit simplex to an error level of $\approx 10^{-14}$ (see \Cref{tab:self_eval}).}}\label{fig:closeeval}
\end{figure}
The quadrature rule~\eqref{eq:Vk_bdry_quad} can lead to a highly-accurate quadrature rule for the representation~\eqref{eq:poincare_Ik} of $\widehat{\mathcal{V}}_k[\rho](\vv{\zeta}_0)$, provided an appropriate distribution of quadrature nodes $\vv{\zeta}_j$ on the boundary $\partial\widehat{\mathcal{C}}$ are selected in the quadrature rule~\eqref{eq:Vk_bdry_quad} for~\eqref{eq:poincare_Ik}. However, simply using a uniform distribution of boundary intervals (e.g.\ for a triangle $\widehat{\mathcal{T}}_0$ parametrized by $Z = Z(t)$, $Z: [0, 2\pi] \to \partial\widehat{\mathcal{C}}$, using interval endpoints equi-spaced in $t$ in each of $[0, 2\pi/3]$, $[2\pi/3, 4\pi/3]$, and $[4\pi/3, 2\pi]$) will not always lead to a high-order approximation of the desired volume potential. Possible loss of accuracy can be understood by observing, in the test depicted in \Cref{fig:closeeval}, that we are building a quadrature rule on $\widehat{\mathcal{C}}$ centered around the target $\vv{\zeta}_*$, with the value at each quadrature node $\vv{\zeta}_j$ in~\eqref{eq:poincare_Ik} being given by the integral $I^k[\rho](\vv{\zeta}, \vv{\zeta}_0)$ along a ray from $\vv{\zeta}_*$ to $\vv{\zeta}_j$. However, an equi-arclength interval distribution does not uniformly cover the angular variable in the coordinate system with origin $\vv{\zeta}_*$, in which context steep gradients arise in $I^k[\rho](\vv{\zeta}, \vv{\zeta}_0)$ with respect to $\vv{\zeta}$, and poor accuracy in the quadrature rule~\eqref{eq:Vk_bdry_quad} can result---a deficiency that we now remedy.

Since potential loss of accuracy in the quadrature rule~\eqref{eq:Vk_bdry_quad} arises due to the described inadequate coverage of the angular variable in the coordinate system with origin $\vv{\zeta}_*$, with this effect being most pronounced at the closest point $\vv{\zeta}_1 \in \partial \widehat{\mathcal{C}}$ to $\vv{\zeta}_*$, we refine the boundary discretization in a vicinity of $\vv{\zeta}_1$ and thus recover volume quadrature nodes that are suitably equi-distributed. After first identifying the parametric location $t_1$ of $\vv{\zeta}_1$ in the parametrization of $\widehat{\mathcal{C}}$, the method proceeds by placing a graded sequence of intervals near to $t_1$. Letting $t_1^-$ and $t_1^+$ denote the parametric location of the nearest (still uniformly-distributed) composite interval endpoints that surround $t_1$, the method introduces additional composite regions with endpoints in the set
\[
    \Pi(t_1) := \left\{t^j: t^j := t_1 \pm |t_1^+ - t_1^-| R^j, \quad j = 1, \ldots, N\right\}; \quad \mbox{with}\quad R = \frac{1}{4}, \quad N = 5.
\]
A depiction of the boundary discretization and associated volumetric quadrature nodes can be seen in \Cref{fig:closeeval}.

\begin{table}[!ht]
    \centering
    \bgroup
    \def\arraystretch{1.2}
    \begin{tabular}{|| c|c|c ||}
        \hline
        $p$ & $\left\|e\right\|_{\ell^\infty}$ & $\widehat{d}_\mathrm{min}$ \\
        \hline
        $4$ & $4.0$e$-15$ & $7.1$e$-02$ \\
        $5$ & $5.0$e$-15$ & $4.8$e$-02$ \\
        $6$ & $7.0$e$-15$ & $3.1$e$-02$ \\
        $7$ & $1.0$e$-14$ & $2.3$e$-02$ \\
        $8$ & $1.1$e$-14$ & $1.9$e$-02$ \\
        $9$ & $1.1$e$-14$ & $1.5$e$-02$ \\
        $10$ & $1.2$e$-14$ & $1.3$e$-02$ \\
        \hline
    \end{tabular}
    \egroup
    \caption{{\em\small Single-cell test of the proposed singular quadrature scheme: evaluation of a volume potential over a straight (mapped) triangle $\Omega$ with vertices located at $(-0.618, -0.312)$, $(-0.825, -0.311)$, and $(-0.802, -0.516)$. For each integer $p$, $\left\|e\right\|_{\ell^\infty}$ denotes maximum error in the volume potential $\mathcal{V}[K_{nm} \circ \vv{R}^{-1}](\vv{r}_0)$ of the $p(p+1)/2$-numbered polynomials $K_{nm} \in \{K_{nm}: 0 \le m \le n, 0 \le n < p\}$, each over all of the $p(p+1)/2$-numbered points $\vv{r}_0 = \vv{R}(\vv{\zeta}_0)$, $\vv{\zeta}_0 \in I_{\mathcal{T}, p}$ (see~\eqref{eq:VR_interp_nodes}). Ground `truth' is a highly-adaptive (and highly-expensive) multidimensional quadrature rule unrelated to the methods of this article, whose associated error is no greater than $\approx 10^{-14}$. The minimum distance from any point $\vv{\zeta}_0$ to $\partial\widehat{\mathcal{T}}_0$ is denoted by $\widehat{d}_\mathrm{min}$.}} \label{tab:self_eval}
\end{table}
\Cref{tab:self_eval} displays results of a test which demonstrates that this refinement strategy for the volume potential yields accuracy of approximately $13$ digits for targets close to cell boundaries. The test consists of evaluation of a Laplace volume potential, $G(\vv{r}, \vv{r}_0) = -\frac{1}{2\pi}\log|\vv{r} - \vv{r}_0|$, over a single (mapped) triangle $\Omega$ at specific points corresponding to Koornwinder interpolation nodes, with associated map denoted by $\vv{R}: \widehat{\mathcal{T}}_0 \to \Omega$. This experiment has direct relevance to evaluation of the volume potential problem~\eqref{eq:volPot} since, as detailed in \Cref{sing_precomp_quad_nodes} below, values of $\mathcal{V}[f]$ at these nodes can be used to produce high-order accurate interpolation of the volume potential throughout a given cell $\mathcal{T}_k$. The results of the experiment summarized in \Cref{fig:Poisson_polydisp} provide a somewhat more challenging test as it includes some target points that lay up to two orders of magnitude closer to the cell boundary.

\subsubsection{Final rule for local corrections}
Having developed singular and near-singular quadrature rules for the one-dimensional integrals $I^k$~\eqref{eq:Ik_def}, we provide the full singular and near-singular rules for the volume potential $\mathcal{V}_k[\rho]$. In the singular case ($\vv{\zeta}_0 \in \mathcal{C}$), inserting the quadrature rule \Cref{eq:Ik_vp_singular} into~\eqref{eq:Vk_bdry_quad} we have
\begin{equation}\label{eq:vp_singular}
\begin{split}
    \widehat{\mathcal{V}}_k[\rho]\left(\vv{\zeta}_0\right) \approx \sum_{j=1}^{PM_{\partial\mathcal{C}}} \left( \sum_{i = 1}^{m_P} t_i \right.&\, G\left(\vv{R}^k(\vv{\zeta}_*+t_i \left(\vv{\zeta}_j-\vv{\zeta}_*\right)),\vv{R}^k(\vv{\zeta}_0)\right)\times \\
    &\left. \vphantom{\sum_{i = 1}^{m_P}} \times J^k( \vv{\zeta}_*+t_i \left(\vv{\zeta}_j-\vv{\zeta}_*\right))\rho\left(\vv{\zeta}_*+t_i\left(\vv{\zeta}_j-\vv{\zeta}_*\right)\right)v_i  \right) \left((\vv{\zeta}_j - \vv{\zeta}_*)\times\vv{\tau}_j\right) w_j.
\end{split}
\end{equation}
Similarly, for the near-singular case ($\vv{\zeta}_0\not\in \widehat{\mathcal{C}}$), inserting the quadrature rule \Cref{eq:Ik_vp_nearsingular} into~\eqref{eq:Vk_bdry_quad} we have
\begin{equation}\label{eq:vp_near_singular}
\begin{split}
    \widehat{\mathcal{V}}_k[\rho]\left(\vv{\zeta}_0\right) \approx \sum_{j=1}^{PM_{\partial\mathcal{C}}} \left( \sum_{i = 1}^{m'_P} \tilde{t}_i \right.&\, G\left(\vv{R}^k\left(\vv{\zeta}_*+\tilde{t}_i\left(\vv{\zeta}_j-\vv{\zeta}_*\right)\right),\vv{R}^k\left(\vv{\zeta}_{0}\right)\right)\times\\
    &\left. \vphantom{\sum_{\ell=1}^{m'_p}} \times J^k(\vv{\zeta}_*+\tilde{t}_i\left(\vv{\zeta}_j-\vv{\zeta}_*\right)) \rho\left(\vv{\zeta}_*+\tilde{t}_i \left(\vv{\zeta}_j-\vv{\zeta}_*\right)\right)\tilde{v}_i \right) \left((\vv{\zeta}_j - \vv{\zeta}_*)\times\vv{\tau}_j\right) w_j.
    \end{split}
\end{equation}

\section{Efficient generation and application of singular corrections}\label{sec:sing_corr_koornwinder}
This section completes the description of a singular quadrature-corrected scheme for the volume potential with an emphasis on \emph{efficiency} that has not as-yet been considered per se: while on the one hand for cells $\mathcal{T}_k$ that are well-separated from the target point $\vv{r}_0$ the smooth quadrature scheme described in \Cref{sec:sparse_func_approx} that leads to the approximation~\eqref{eq:Tk_smooth_quad} is accurate (and amenable to FMM acceleration), and the resulting rules for singular and near-singular quadrature corrections (respectively~\eqref{eq:vp_singular} and~\eqref{eq:vp_near_singular}) from \Cref{sec:sing_quad} accurately evaluate contributions to $\mathcal{V}[f]$ from nearby cells, on the other hand we have still yet to describe an efficient scheme. By `efficient' we really mean simultaneously that
\begin{enumerate}[(i)]
    \item The method is `node-efficient', that is, efficient with respect to the required number of degrees of freedom per cell,
    \item The local corrections are cheap to \emph{apply} (i.e.\ they are cheap relative to the FMM call), and
    \item The local corrections are cheap to \emph{generate}.
\end{enumerate}
To motivate our approach we note that simply applying the methods of \Cref{sec:sparse_func_approx} to the source density $f$ is \emph{inefficient} since, while accurate, the associated points where the source density $f$ are required, the singular quadrature nodes, are both numerous and dependent on the target itself. Our approach, rather, is to apply the singular quadrature methods in \Cref{sec:sing_quad} to the Koornwinder basis elements (which, recalling \Cref{sec:triangles} provide a high-order basis for approximation of arbitrary smooth functions on $\mathcal{T}_k$). The singular corrections are highly efficient for repeated application of the Newton potential since they are (a) \emph{Local} in the sense that the only cells which contribute to a singular or near-singular correction at a target $\vv{r}_0$ are those cells which are a subset of $\mathcal{C}^\mathrm{self}(\vv{r}_0) \cup \mathcal{C}^\mathrm{near}(\vv{r}_0)$ and are also (b) \emph{Data-sparse}, as a linear map from source function $f$ values at Koornwinder interpolation nodes directly to local correction values.

In more detail, as a consequence of the local correction methods described in this section and the smooth oversampled quadratures described in \Cref{sec:boxes} we conclude that the total degrees of freedom (source function evaluation points) of the scheme for~\eqref{eq:volPot} number
\begin{equation}\label{eq:ndofs}
    \mathrm{ndofs} = \frac{p(p+1)}{2} N_t + p^2 N_b,
\end{equation}
while the source points for the FMM number
\begin{equation}\label{eq:nsrcs}
    \mathrm{nsrcs} = \frac{q(q+1)}{2}N_t + q^2 N_b, \quad q \ge p,
\end{equation}
both for the order $p$ scheme over a mesh with $N_t$ triangular regions and $N_b$ boxes. Denoting by $\Pi_k(\vv{r}_0): \mathbb{R}^{N_s} \to \mathbb{R}$ the linear map from $N_s$ source function values on a cell $\mathcal{C}_k$ to scalar correction (i.e.\ the high-order accurate value of $\mathcal{V}_k$ obtained via singular quadrature less the inaccurate contribution from the oversampled smooth quadrature described in \Cref{sec:triangles} for triangles and \Cref{sec:boxes} below for boxes) for the volume potential on $\mathcal{C}_k$ evaluated at $\vv{r}_0$, we have either $N_s = p(p+1)/2$ (in the case $\widehat{\mathcal{C}}_k = \widehat{\mathcal{T}}_0$) or $N_s = p^2$ (in the case $\widehat{\mathcal{C}}_k = \widehat{\mathcal{B}}_0$) with, of course, an $\mathcal{O}(1)$ number of nontrivial correction maps $\Pi_k(\vv{r}_0)$ per target point due to the $\mathcal{O}(1)$ cardinality of $\mathcal{C}^\mathrm{near}_k(\vv{r}_0)$ (see \Cref{sec:meshing}). Assembling all nontrivial corrections $\Pi_k$ into a sparse matrix results in a final singular correction scheme where the overwhelming majority of computational effort is spent in the highly-efficient FMM stage, and which requires a conservative number of source evaluation nodes; the overall performance of the solver is demonstrated in \Cref{tab:fmm_timing}.
\begin{remark}
    Numerical experiments show that a fixed, twofold ($q=2p$) oversampling allows for smooth quadrature error for $\mathcal{V}_k[f]$ that is, as desired, dominated by interpolation error of $f$ for the smooth (Stokes and Laplace) kernels and for the error levels presented in this article (thus, for such kernels the selection $q = 2p$ was made in producing the numerical results), though see also~\cite{Greengard:21} for other oversampling schemes. In particular, to achieve higher accuracies than those presented in our numerical experiments, an increase in the fixed oversampling rate is sometimes useful---increasing slightly the number of FMM source points, but not the overall number of degrees of freedom, per~\eqref{eq:ndofs} and~\eqref{eq:nsrcs}. Furthermore, oversampling and adequate identification of an appropriately-sized near-field region is more critical to obtain high accuracies for more sharply-peaked kernels, such as Helmholtz kernels. We leave detailed consideration of these matters to future study.
\end{remark}

Point (iii) above relates to a question of perennial concern in potential theoretic methods: the efficient generation of corrections for singular and near-singular targets (see reference~\cite{Greengard:21} where adaptive quadrature is performed for each near-singular target, at significant cost, and see also references~\cite{Bremer:12, Bremer:13, BrunoGarza:20, Klockner:20, Perez:19, Faria:21}); the cost of generating such corrections can be burdensome, even while the resulting per-use costs, say, in an iterative solver seem highly favorable. (In our context the cost is due both to the generation of target-adapted quadrature nodes and weights and also to the evaluation of the orthogonal polynomial family at these nodes.) We address this issue by showing that as a result of our use of mappings from a common reference-space domain, the quadratures and values of orthogonal polynomial values in reference-space can be re-used across elements even as the mappings vary, with the favorable implication that the only required per-element computation for generating the (near-)singular corrections consist of evaluation of the Green function.

\subsection{Triangles and Koornwinder systems}
For a given target $\vv{r}_0$ and for each nearby cell $\mathcal{T}_k$, i.e.\ $\mathcal{T}_k$ satisfying $\mathcal{T}_k \subset \mathcal{C}^\mathrm{self}(\vv{r}_0) \cup \mathcal{C}^\mathrm{near}(\vv{r}_0)$, the method begins by approximating the source function $f$ on $\mathcal{T}_k$ using the $p$-th order Koornwinder expansion
\begin{equation}\label{eq:f_koornwinder_approx}
    f(\vv{R}^k(\xi, \eta)) \approx \sum_{n=0}^{p-1}\sum_{m=0}^n a_{nm}^k K_{nm}(\xi, \eta),
\end{equation}
where the coefficients $\vv{a}^k$ solve the system
\begin{equation}\label{eq:amap_koorn}
    \vv{V}_p\vv{a}^k = \vv{f}^k_p,
\end{equation}
with $(\vv{f}_p^k)_{i} = f(\vv{R}^k(\xi_{p,i}, \eta_{p, i})),\, i = 1, 2, \ldots, N_p$, the vector of function values at the Koornwinder interpolation nodes of~\eqref{eq:VR_interp_nodes}. Substituting this expansion into the volume potential~\eqref{eq:Vk_hatVk_equiv} we obtain the high-order approximation to $\mathcal{V}_k[f](\vv{r}_0)$,
\begin{equation}\label{eq:maptri_singular_quad}
    \begin{split}
        \mathcal{V}_k[f](\vv{r}_0) &\approx \int_{\widehat{\mathcal{T}}_0} G(\vv{R}^k(\xi, \eta), \vv{R}^k(\xi_0^k, \eta_0^k)) \sum_{n=0}^{p-1} \sum_{m=0}^n a_{nm}^k K_{nm}(\xi, \eta) J^k(\xi, \eta) \,\d\xi\d\eta\\
        &= \sum_{n=0}^{p-1} \sum_{m=0}^n a_{nm}^k \int_{\widehat{\mathcal{T}}_0} G(\vv{R}^k(\xi, \eta), \vv{R}^k(\xi_0^k, \eta_0^k) K_{nm}(\xi, \eta) J^k(\xi, \eta) \,\d\xi\d\eta\\
        &= \sum_{n=0}^{p-1} \sum_{m=0}^n a_{nm}^k \widehat{\mathcal{V}}_k[K_{nm}](\vv{\zeta}_0^k),
    \end{split}
\end{equation}
where
\[
    \vv{\zeta}_0^k = \left( \xi_0^k, \eta_0^k\right)^T
    = \left(\vv{R}^k\right)^{-1}(\vv{r}_0)
\]
denotes the $\mathcal{T}_k$-reference-space location of the target $\vv{r}_0$. Writing the volume potential $\mathcal{V}_k[f](\vv{r}_0)$ in terms of reference-space potentials $\widehat{\mathcal{V}}_k$ allow for singular and near-singular correction of the smooth quadrature rule~\eqref{eq:Tk_smooth_quad} for cells $\mathcal{T}_k$ that lay close to $\vv{r}_0$. The corrections can be pre-computed for each element in the $p$-th order Koornwinder using the methods of \Cref{sec:sing_quad}---specifically, if $\vv{r}_0 \in \mathcal{T}_k$ then rule~\eqref{eq:vp_singular} is used to compute $\widehat{\mathcal{V}}_k[K_{nm}](\vv{\zeta}_0^k)$ while otherwise rule~\eqref{eq:vp_near_singular} is used.

This completes the description of how singular and near-singular corrections can be generically pre-computed; in what follows we consider optimizations that can be obtained for the fortunately-typical case of fixed reference-space target locations. We describe here the case for prescribed reference-space target points $\vv{\zeta}_0 = (\xi_0, \eta_0) \in I_{\mathcal{T}, p}$, where $I_{\mathcal{T}, p}$ are the Vioreanu-Rokhlin interpolation nodes (see also \Cref{sec:triangles}); this ensures the possibility of high-quality interpolation of the resulting function $\mathcal{V}[f]$ throughout the domain (which is useful, for example, in the context of a non-linear or time-dependent PDE)---see also the experiment in \Cref{sec:results} corresponding to \Cref{tab:fmm_timing}. The ideas are not restricted to a specific set of interpolation nodes, and the method proposed in \Cref{sing_precomp_quad_nodes} could be used e.g.\ for fixed reference-space locations on the curved boundary component of a curvilinear cell---i.e.\ for evaluation of $\mathcal{V}[f](\vv{r}_0)$ with $\vv{r}_0 \in \partial \Omega$. Optimized pre-computation of singular corrections are generalized later to near-singular targets.

\subsubsection{Efficient pre-computation of singular and near-singular corrections}\label{sing_precomp_quad_nodes}
An observation that leads to substantial efficiency gains for the proposed method is that the quadratures for \emph{every} cell $\mathcal{T}_k$ occur on the same reference cell $\widehat{\mathcal{T}}_0$. Since every triangular mesh cell is mapped to the standard simplex $\widehat{\mathcal{T}}_0$ and the target nodes are fixed in reference space, the same (expensive to generate) singular quadratures and Koornwinder polynomial values at these nodes can be used for each target point \emph{for every cell}. To see the implications for computational cost savings, observe that from~\eqref{eq:maptri_singular_quad},
\[
    \mathcal{V}_k[f](\vv{r}_0) \approx \sum_{n=0}^{p-1} \sum_{m=0}^n a_{nm}^k \widehat{\mathcal{V}}_k[K_{nm}](\vv{\zeta}_0),
\]
with
\begin{equation}\label{eq:Vnm_k_def}
  \widehat{\mathcal{V}}_k[K_{nm}](\vv{\zeta}_0) = \int_{\widehat{\mathcal{T}}_0} G(\vv{R}^k(\vv{\zeta}), \vv{R}^k(\vv{\zeta}_0)) K_{nm}(\vv{\zeta}) J^k(\vv{\zeta}) \,\d A(\vv{\zeta}).
\end{equation}
Clearly, the singular point is $\vv{\zeta} = \vv{\zeta}_0$ \emph{for every cell} $\mathcal{T}_k$, meaning that the singular quadrature weights and nodes are also identical even as the mapping $\vv{R}^k$ and associated Jacobian $J^k$ change. Applying the singular quadrature rule~\eqref{eq:vp_singular} with $\rho = K_{nm} J^k$ and $\widehat{\mathcal{C}} = \widehat{\mathcal{T}}_0$ results in a set $\lbrace(\vv{\chi}_i, \omega_i): i = 1, \ldots, N^Q\rbrace$, $N^Q = PM_{\partial S}m_P$, of pairs of quadrature nodes $\vv{\chi}_i = (\xi_i, \eta_i)$ and weights $\omega_i$ pairs adequate for discretization of $\widehat{\mathcal{V}}_k$ that are independent of $k$ (here, for appropriate $j$ and $\ell$ in~\eqref{eq:vp_singular} the weights $\omega_i$ are given by $\omega_i = t_\ell v_\ell((\vv{\zeta}_j - \vv{\zeta}_*)\times \vv{\tau}_j)$). Re-writing~\eqref{eq:vp_singular} with this notation we have the singular corrections given by the $\omega_j$-weighted inner product of Green function and Jacobian values with Koornwinder polynomial values,
\begin{equation}\label{eq:singcorr_inner_prod}
  \widehat{\mathcal{V}}_k[K_{nm}](\vv{\zeta}_0) \approx \sum_{i=1}^{N^Q} G(\vv{R}^k(\vv{\chi}_i), \vv{R}^k(\vv{\zeta}_0)) J^k(\vv{\chi}_i) K_{nm}(\vv{\chi}_i) \omega_i,
\end{equation}
from which it is clear that the only quantities that need to be recomputed for each cell are the mapped Green function values $G(\vv{R}^k(\vv{\chi}_i), \vv{R}^k(\vv{\zeta}_0)) J^k(\vv{\chi}_i)$; note, further, that $J^k$ is constant for straight triangles. It is \emph{essential} in the independence of the node-weight pairs $(\vv{\chi}_i, \omega_i)$ with respect to the mapping $k$ (and thereby to the cell) that for a given target $\vv{\zeta}_0$ the singular quadrature rule developed in \Cref{sec:sing_quad} is \emph{determined solely by the small-argument asymptotic behavior of the kernel function}.

The additional challenge that unstructured meshes pose for \emph{near-singular} corrections as opposed to the singular corrections discussed previously is that near-singular target locations $\vv{r}_0 = \vv{R}^j(\vv{\zeta}_0^j)$ (where $\vv{\zeta}_0^j = \vv{\zeta}_0$ is, say, one of the Koornwinder interpolation nodes in cell $\widehat{\mathcal{T}}_j$---\emph{which are always in fixed locations in the $\widehat{\mathcal{T}}_j$-reference space}) unfortunately do lay at \emph{arbitrary} reference space locations $\vv{\zeta}_0^k$ relative to the source cell $\mathcal{T}_k$ ($k \neq j$) under consideration. For this reason it is no longer possible to repeatedly use the same fixed set of reference-space target points to directly generate the corrections as for the $j = k$ case described previously; fortunately, however, similar ideas are still applicable with essentially the same effect. Despite the arbitrary (reference-space) location of near-singular target points, the quadrature rule weights and nodes (and hence the required values of the polynomials $K_{nm}$) are still determined entirely by the resulting star-point $\vv{\zeta}_*$ (which is itself, in turn, determined via \Cref{eq:rstar_def})---see \Cref{fig:nearsingquad_nodes} for a depiction of this dependence. As a result, one can re-use the quadratures generated for a fixed set of star-points laying along the boundary $\partial \widehat{\mathcal{T}}_0$, selecting the star-point $\vv{\zeta}_*$ that lays closest to the solution to \Cref{eq:rstar_def}. Utilizing a fixed-size (independent of the number of cells) list of quadrature nodes, weights, and Koornwinder values $K_{nm}$, the quadrature rules for near-singular targets can be efficiently computed by means of a simple lookup table and, again, the formula~\eqref{eq:singcorr_inner_prod} (where, of course, the mapped Green function still needs to be re-evaluated, as before). The already-modest storage costs of these quadrature rules in such a scheme could be further limited by storing rules corresponding to star points on only one side of $\partial\widehat{\mathcal{T}}_0$.

\subsection{Boxes and tensor-product Chebyshev systems: smooth quadratures and singular corrections}\label{sec:boxes}
In this section, we describe smooth quadrature methods and efficient singular and near-singular corrections of these quadratures for source boxes $\mathcal{C}_{k+N_t} = \mathcal{B}_k$, $k = 1, \ldots, N_b$ (see \Cref{eq:Ck_Tk_Bk}). Every box $\mathcal{B}_k$ is mapped from a single reference box $\widehat{\mathcal{B}}_0 = [-1, 1]^2$, with mappings in the uniform grid context of this article taking the simple form
\begin{equation}\label{eq:Rk_box}
    \vv{R}^k(\vv{\zeta}) = \vv{o}_{k} + h\vv{\zeta}, \quad\mbox{where}\quad \vv{\zeta} = \begin{pmatrix}\xi, \eta\end{pmatrix}^T,\quad\mbox{with} \quad \xi,\eta \in [-1, 1].
\end{equation}

\textbf{Smooth quadrature and interpolation.}
As is standard for such regular regions, we represent the local source distribution $f$ using the truncated tensor-product Chebyshev series expansion
\begin{equation}\label{eq:fk_basis}
    f(\vv{r}) = f(\vv{R}^k(\zeta)) \approx \sum_{n=0}^{p-1}\sum_{m=0}^{p-1} f^k_{nm}T_{nm}(\xi, \eta), \quad\mbox{for}\quad \vv{r} = \vv{R}^k(\vv{\zeta}) \in \mathcal{B}_k,
\end{equation}
where $T_{nm}$ is defined by the tensor-product Chebyshev polynomial
\begin{equation}\label{eq:Tnm}
    T_{nm}(\vv{\zeta}) = T_{nm}(\xi, \eta) = T_n(\xi)T_m(\eta),
\end{equation}
with $T_n(t)$ denoting the Chebyshev polynomial of degree $n$ defined for $-1 \le t \le 1$. We denote by $I_{\mathcal{B}, p}$ the set of interpolation nodes for the truncated series representation~\eqref{eq:fk_basis}, which are defined as tensor-product Chebyshev nodes, i.e.\ roots of the Chebyshev polynomial $T_p$.
It is well-known that for smooth $f$ this is a $p$-th order accurate expansion with series coefficients $f^k_{nm}$ decaying rapidly as either of $n, m$ increase, and, further, that the associated Clenshaw-Curtis smooth quadrature rule is suitable for targets $\vv{r}_0$ well-separated from $\mathcal{B}_k$,
\[
    \int_{\mathcal{B}_k} G(\vv{r},\vv{r}_0)f(\vv{r})\, \mathrm{d}\vv{r} \approx \sum_{I=1}^{p^2} G(\vv{r}_I,\vv{r}_0)f(\vv{r}_I)w_I,
\]
and, further, that is compatible with fast summation methods. Mirroring the upsampling procedure described in \Cref{sec:triangles} the method uses a smooth quadrature of order $q$, $q \ge p$, for an interpolant~\eqref{eq:fk_basis} of order $p$ so that source approximation error is dominant, with the Chebyshev expansion~\eqref{eq:fk_basis} with associated interpolation nodes $I_{\mathcal{B},p}$ upsampled to the interpolation nodes $I_{\mathcal{B},q}$.

\textbf{Singular and near-singular corrections.}
For singular and near-singular evaluation, similar to the case for triangular regions, we construct a linear map from function values $f_{nm}^k$ in \Cref{eq:fk_basis} to corrected quadratures $\mathcal{V}_k(\vv{r}_0)$. The singular corrections are generated using the same Poincare-based singular integration technique used for triangular regions, applied to each element of the smooth basis in~\cref{eq:fk_basis}. Similarly to~\eqref{eq:maptri_singular_quad} we find
\begin{equation}\label{eq:mapbox_singular_quad}
    \mathcal{V}_k[f](\vv{r}_0) \approx \sum_{n=0}^{p-1}\sum_{m=0}^{p-1} f^k_{nm}\widehat{\mathcal{V}}_k[T_{nm}](\vv{\zeta}_0^k),
\end{equation}
where $\widehat{\mathcal{V}}_k[\rho]$ is defined by \Cref{eq:Vkdef} and where from~\eqref{eq:Rk_box} we have $\vv{\zeta}_0^k = (\vv{r}_0 - \vv{o}_k)/h$.
Crucially, by translational invariance $G(\vv{r}, \vv{r}_0) = G(\vv{r} - \vv{r}_0)$ and from the mapping~\eqref{eq:Rk_box} we have
\begin{equation}\label{Vkhat_box}
    \widehat{\mathcal{V}}_k[T_{nm}](\vv{\zeta}_0^k) = \int_{\widehat{\mathcal{B}}_0} G(h(\vv{\zeta} - \vv{\zeta}_0^k)) T_{nm}(\vv{\zeta})h^2 \,\d\vv{\zeta},
\end{equation}
where we emphasize that the expression for $\widehat{\mathcal{V}}_k$ is independent of $k$.

It follows from~\eqref{Vkhat_box} and the fact that the reference-space interpolation and smooth quadrature nodes $I_{\mathcal{B}, p}$ for $\mathcal{B}_k$ are identical for all boxes, that the singular corrections can be computed once and re-used for all boxes as a lookup table. A similar argument leads to lookup tables for near-singular target points that arise at interpolation nodes of nearby boxes: since the grid is structured these need only be computed once (this remains true even with adaptivity, merely the size of the still-finite look-up table expands). Finally, for near-singular targets laying in triangles that lay in close proximity to boxes, corrections can be easily generated with the added efficiency (since, in this context the map for every box is the same) that the values of the Newton potential in the near-field can be re-used across boxes and indeed interpolation to machine precision of the volume potential is possible for generating near-singular corrections.

\section{Numerical results}
\label{sec:results}

This section presents demonstrations of the character of the proposed numerical methods, applied to the Poisson and modified Helmholtz (Yukawa) equations. Together with the Stokes equations (see \Cref{fig:Stokes} for a solution of this equation with the proposed methodology) these equations represent the major classes of constant-coefficient elliptic PDEs to which our methods are applicable, even though, as noted in \Cref{sec:Conclusions}, our methods possess broader applicability. We first show the entirely routine and expected convergence of the associated homogeneous problem, then consider test cases that involve use of the volume potential, demonstrating both convergence and the asymptotic costs of the method with the use of the FMM.

\begin{figure}[!ht]
    \centering
    \includegraphics[height=0.23\textwidth]{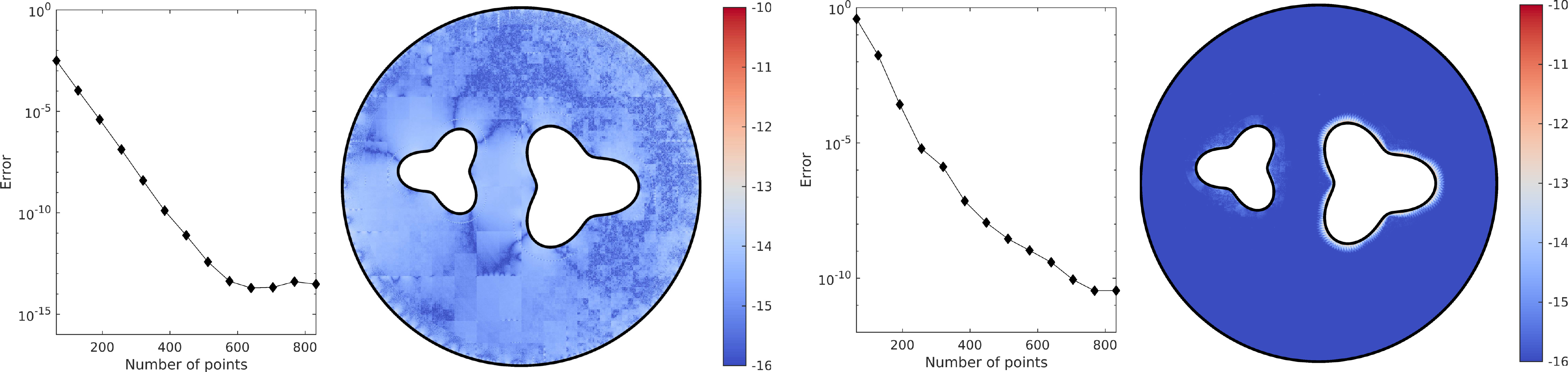}
    \caption{{\em\small Laplace and modified Helmholtz Dirichlet boundary value problem (BVP) test problems, with solution (and hence Dirichlet boundary values) given by point sources (for Helmholtz, $\lambda = 10$) at the center of each inclusion. Left (resp.\ Right): convergence plot of the collocation scheme for the Laplace (resp. Helmholtz) BVP with respect to the total number $N$ of collocation points; error plots correspond to $N=832$.}}
    \label{fig:BIEerror}
\end{figure}

As a preliminary test, we demonstrate the expected convergence of the numerical solver for the integral equations arising in the augmented boundary value problem \Cref{eq:ellipt_augment}. The boundary integral equations are discretized using spectral Kress quadratures, the resulting linear systems are solved using GMRES with a relative residual tolerance of $10^{-15}$, and, for evaluation of the layer potentials spectrally-accurate schemes~\cite{barnett2015spectrally} are used that yield high accuracy in the numerical solution arbitrarily close to the boundary. \Cref{fig:BIEerror} validates the expected exponential convergence of the numerical solution as the total number of quadrature nodes increases. In the remainder of the numerical results, a sufficient number of boundary integral quadrature nodes are selected so that the error from the homogeneous component of the solver does not dominate.

\begin{remark}
Several Poisson examples in this section concern the inhomogeneity (the inhomogeneity being modified in an obvious manner for Helmholtz problems) given by
\begin{equation}\label{eq:full_domain_mfg_soln}
        f(x, y) = 6\sin(6x) + 8\cos(8(y + \frac{1}{10})) + 4(x^2 + y^2)\sin(4xy) + 3\cos(3x)\sin(3y)
\end{equation}
with associated solution
\begin{equation}\label{eq:full_domain:mfg_soln_u}
        u(x, y) = \frac{1}{6} \sin(6x) + \frac{1}{8} \cos(8(y + \frac{1}{10})) + \frac{1}{4} \sin(4xy) + \frac{1}{6}\cos(3x)\sin(3y).
\end{equation}
We emphasize that while $f(x, y)$ is defined and is smooth for all $(x, y)$, only function values at interior points $(x, y) \in \Omega$ are used in the solution process.
\end{remark}

\begin{figure}[!ht]
    \centering
    \includegraphics[width=\textwidth]{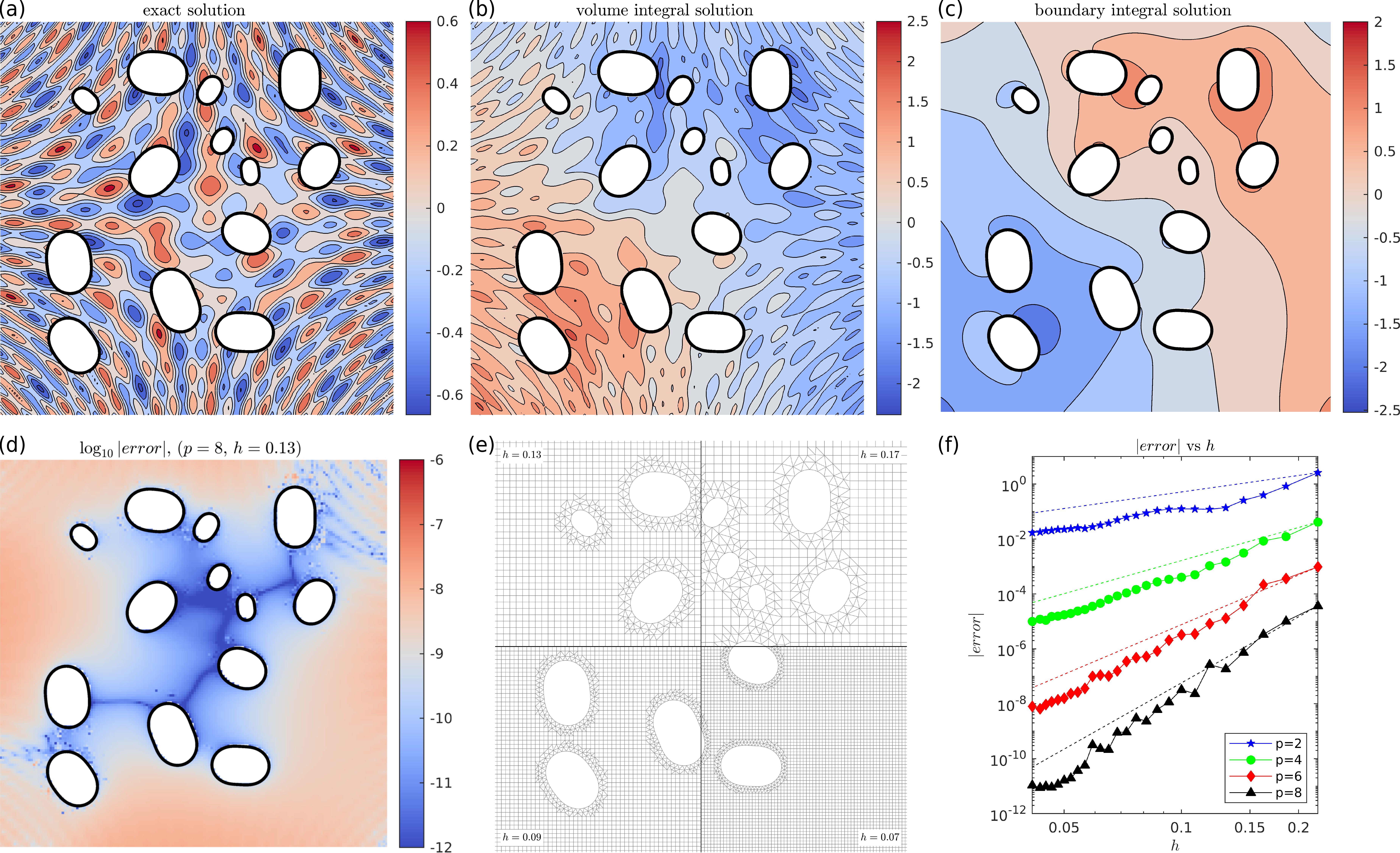}
    \caption{{\em\small Poisson problem simulation; solution and errors displayed in the subregion $[-4.2, 4.2]^2$ of a polydisperse domain. The three top subplots~(a-c) depict, respectively, the contour lines for the associated solution $u$ to \Cref{eq:ellipt}, the particular solution $u_P$ arising from \Cref{eq:volPot}, and the solution $u_H$ to the boundary value problem~\eqref{eq:ellipt_augment}. Subplot~(d) is a typical plot of the base-$10$ logarithm of the error at points in a uniform evaluation grid, here corresponding to the solver run with $p=8$ and $h=0.13$. Subplot~(e) displays four computational meshes at various levels of $h-$refinement. Subplot~(f) demonstrates the convergence (with respect to gridsize $h$) of the numerical solution produced by the order $p=2,4,6$ and $8$ versions of the scheme; errors are measured on the same $100\times 100$ uniform target grid over $[-6.5, 6.5]^2$ in all simulations, and dashed lines depict the expected order of convergence.}}
    \label{fig:Poisson_polydisp}
\end{figure}

We consider a Dirichlet problem for the Poisson equation in the presence of a polydisperse system $\Omega$ of inclusions, in the interior of a circle of radius $6.5$ units. The solution and its numerical components are displayed in the region $[-4.2, 4.2]^2$ in \Cref{fig:Poisson_polydisp} (the solution in all of $x^2 + y^2 < 6.5^2$ is not displayed, but the error and therefore the convergence plot in panel (f) is computed over the entire solution domain). The Poisson solution is tested on a $100 \times 100$ uniform grid of target points that lay inside $\Omega$; the maximum error is plotted in panel (f) of \Cref{fig:Poisson_polydisp}. This experimental setup can be seen as a somewhat more challenging test case for the volume potential scheme since it results in targets that can lay arbitrarily close to cell boundaries. Nevertheless we see convergence rates consistent with expectations from approximation theory down to a level on the order of $10^{-12}$. In \Cref{fig:Poisson_polydisp} panel (d) we observe that the error grows as $x$ and/or $y$ increase which is expected in view of the locations of sharper gradients of $f(x, y)$. \Cref{tab:Poisson_polydisp} assembles relevant mesh data for this problem as the discretization is refined and demonstrates the expected linear (resp.\ quadratic) growth of the number $N_t$ (resp.\ $N_b$) of boundary fitted cells (resp.\ regular boxes).

\begin{table}[!ht]
    \centering
    \begin{tabular}{|| c|c|c|c|c|c|c|c|c|c ||}
        \hline
        $N$ & $70$ & $100$ & $130$ & $160$ & $190$ & $220$ & $250$ & $280$ & $310$\\
        \hline
        $h$ & $0.1857$  &  $0.1300$  &  $0.1000$  &  $0.0813$  &  $0.0684$  &  $0.0591$  &  $0.0520$  &  $0.0464$  &  $0.0419$\\
        $N_t$ & $2329$  &  $3293$  &  $4230$  &  $5205$  &  $5698$  &  $6643$  &  $8049$  &  $8531$  &  $8965$\\
        $N_b$ & $2476$  &  $5774$  &  $10408$  &  $16358$  &  $23618$  &  $32227$  &  $42162$  &  $53403$  &  $65932$ \\
        $hN_t / D_1$ & --- & $0.99$ & $0.98$ & $0.98$ & $0.90$ & $0.90$ & $0.97$ & $0.92$ & $0.87$ \\
        $h^2N_b / D_2$ & --- & $1.14$ & $1.22$ & $1.27$ & $1.29$ & $1.32$ & $1.34$ & $1.35$ & $1.36$ \\
        \hline
    \end{tabular}
    \caption{{\em\small Total number of triangular cells $N_t$ and regular boxes $N_b$ at various levels of $h-$refinement ($h = 13/N$) for the polydisperse test problem depicted in \Cref{fig:Poisson_polydisp}: the interior of a circle of radius $6.5$ units and exterior to the elliptical inclusions. The last two rows demonstrate the linear and quadratic growth of $N_t$ and $N_b$, respectively, as $h \to 0$; the ratios $D_1$ and $D_2$ are normalization constants $D_1 = 0.1857 \cdot 2329$ and $D_2 = 0.1857^2 \cdot 2476$.}}\label{tab:Poisson_polydisp}
\end{table}

\begin{figure}
    \centering
    \includegraphics[width=0.40\textwidth]{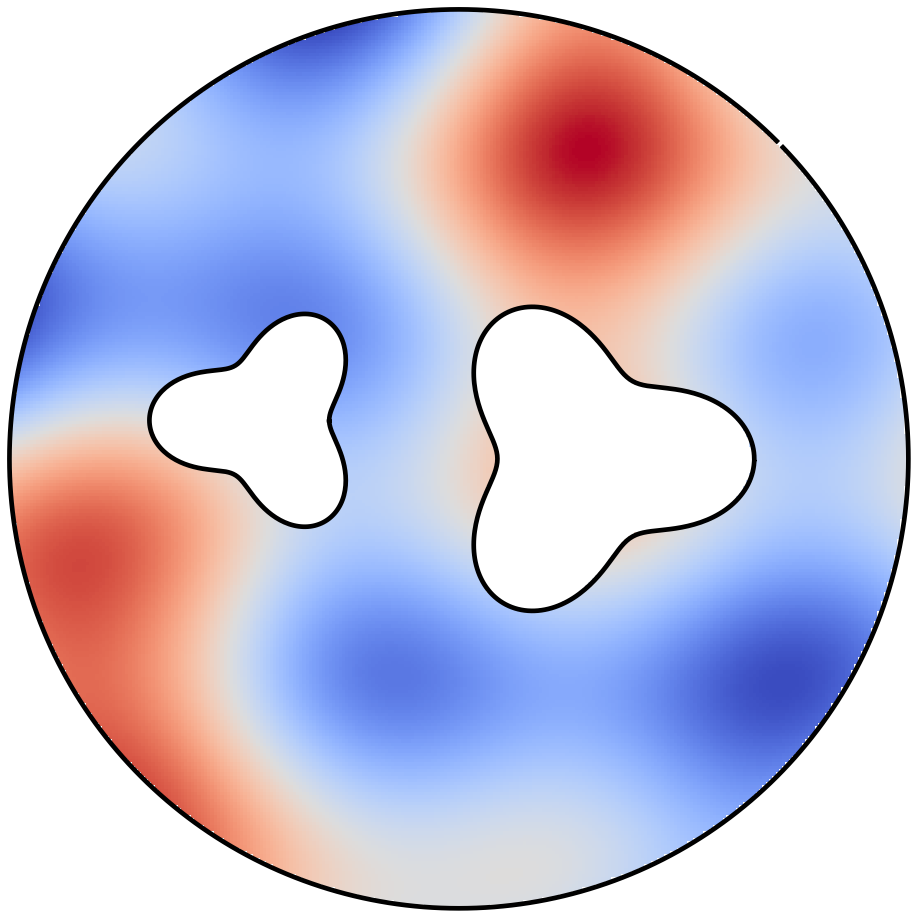}
    \includegraphics[width=0.48\textwidth]{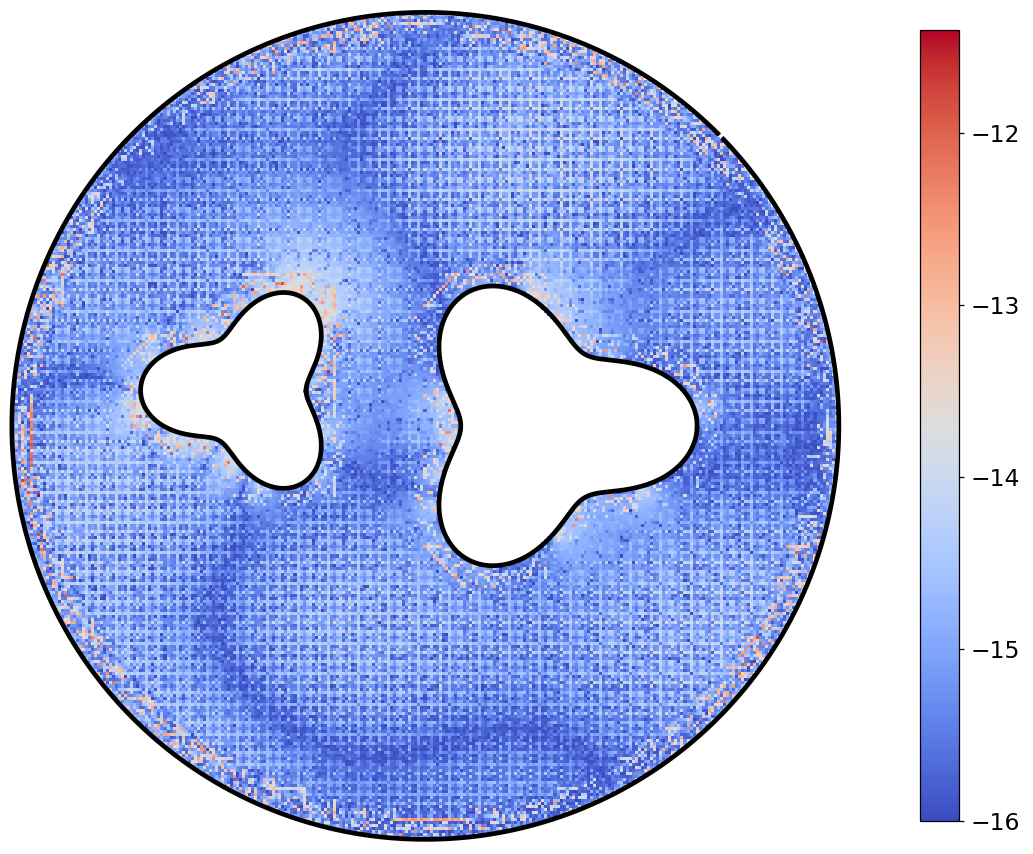}
    \caption{{\em\small Modified Helmholtz ($\lambda^2 = 100$) example solution with maximum value $0.56$ (left) and error with color bar indicating the base-$10$ logarithm of the error (right); solver parameters used were $p = 10$, $h = 2.2/N$ with $N = 90$ (resulting in $N_b = 3633$ and $N_t = 1994$). The maximum error over approximately $51000$ uniformly-spaced target points is $\left\|\vv{e}\right\|_{\ell^\infty} = 3.9 \cdot 10^{-12}$, while the discrete $L^2$ error is $\left\|\vv{e}\right\|_{\ell^2} = 6.9 \cdot 10^{-14}$.}}
    \label{fig:modhelm}
\end{figure}

It is also of interest to solve the inhomogeneous modified Helmholtz equation (cf.\ \Cref{sec:background} with $L = -\Delta + \lambda^2$), a problem with diverse applications, e.g.\ in wave scattering and elliptic time marching. The Green function of this operator shares the same singular kernel behavior as the classical Newton potential, and so the same singular and near-singular quadrature rules used there apply. In \Cref{fig:modhelm} we show the error in the proposed method applied to the modified Helmholtz equation with manufactured solution given in~\eqref{eq:full_domain:mfg_soln_u} over a region bounded by a circle of radius $1$. The method presented is quite successful for low to moderate real values of $\lambda$ with no modifications. However, both the accuracy of smooth quadrature rules and efficacy of the FMM are strained as $\lambda$ is increased; a more complete study of the optimal relation between $\lambda$, the smooth quadrature rule and associated upsampling rate $q$, and the near-field selection region will be presented at a later date.

We turn next to performance demonstrations of the proposed method, which are summarized in \Cref{tab:fmm_timing}. This experiment tests the solution of the Dirichlet Poisson problem in the region depicted in \Cref{fig:BIEerror} (which features an enclosing circle of radius $3.5$ units), with source function and associated solution given by \Cref{eq:full_domain_mfg_soln} and for volume target points at all volume interpolation nodes (i.e.\ for each cell $\mathcal{C}_k$ the appropriate (mapped) points in either of $I_{\mathcal{T}, p}$ and $I_{\mathcal{B}, p}$), and it allows us to make several observations and conclusions about the character of the solver. These experiments serve first to demonstrate that the costs associated with the sparse correction step of the method are, as desired, negligible in comparison to the FMM step of the method, which can be verified with the data in columns marked ``$\%$fmm''. The timings clearly confirm the expected linear scaling of the method. The experiments also serve to demonstrate the costs and error quantities associated with producing the solution at all interpolation nodes in the domain decomposition scheme; the values of the solution at such points suffice to efficiently and accurately produce volume potential values at arbitrary targets $\vv{r}_0 \in \Omega$, but more crucially they can be easily seen as the values required for the solution of various nonlinear and time-dependent PDEs, by e.g.\ iteration or time-stepping. The error $\left\|e_1\right\|_{\ell^\infty}$ in the second-to-last column of \Cref{tab:fmm_timing} is the maximum error in the solution $u$ of the problem~\eqref{eq:ellipt} evaluated across all ndof volume quadrature points and is consistent with the expected $5$-th order convergence for the $p=5$ scheme. Finally, the experiment demonstrates the expected linear (resp.\ quadratic) dependence of the number $N_t$ (resp.\ $N_b$) of triangular regions (resp.\ boxes) on the mesh gridsize $h$, which for this experiment is given by $h = 9/N$. Results are similar for the scheme with other orders, and results are not included for the sake of brevity. All timings were produced using a Python implementation on a single core (i.e.\ with multithreading explicitly disabled) of an Exherbo Linux workstation with an AMD Ryzen-5 5600X CPU\@. A detailed presentation of performant generation of singular and near-singular corrections as detailed in \Cref{sing_precomp_quad_nodes} will be reported at a later date.

It is also of interest to obtain the solution at arbitrary points in the domain via interpolation, a straightforward task when (as in the context of the test of the previous paragraph) evaluation points are precisely the interpolation nodes $I_{\mathcal{T},p}$ and $I_{\mathcal{B},p}$. The error $\left\|\vv{e}_2\right\|_{\ell^2}$ in the solution thus obtained by Koornwinder and Chebyshev interpolation is listed in the last column of \Cref{tab:fmm_timing}, where $\left\|\vv{e}\right\|_{\ell^2} = \left( \sum_j h^2 (\vv{e})_j^2\right)^{1/2}$ denotes the discrete $L^2$ norm of the error $\vv{e}$, and where $h = 0.045$ is the gridsize of the uniform evaluation grid (with $\approx 17,000$ points in $\Omega$) onto which solution values are interpolated.
\begin{table}
    \centering
    \begingroup
    \setlength{\tabcolsep}{4pt}
    \begin{tabular}{|| c|c|c|c|c|c|c|c|c|c|c ||}
        \hline
        $N$ & $N_t$ & $N_b$ & $T^\mathrm{fmm}_{\Gamma}$ & $T_\Gamma^\mathrm{corr}$ & $\%$fmm & $T_\Omega^\mathrm{fmm}$ & $T_\Omega^\mathrm{corr}$ & $\%$fmm & $\left\|\vv{e}_1\right\|_{\ell^\infty}$ & $\left\|\vv{e}_2\right\|_{\ell^2}$\\
        \hline
        $40$ & $828$ & $286$ & $1.79$e-01 & $1.1$e-04 & $99.9$ & $4.30$e-01 & $2.44$e-03 & $99.4$ & $8.2$e-05 & $1.3$e-03 \\
        $60$ & $1147$ & $946$ & $3.05$e-01 & $1.1$e-04 & $99.9$ & $9.54$e-01 & $5.07$e-03 & $99.5$ & $1.1$e-05 & $7.1$e-05 \\
        $80$ & $1440$ & $1922$ & $6.44$e-01 & $1.4$e-04 & $99.9$ & $1.53$e+00 & $7.61$e-03 & $99.5$ & $2.6$e-06 & $1.3$e-05 \\
        $100$ & $1858$ & $3227$ & $8.25$e-01 & $1.1$e-04 & $99.9$ & $2.76$e+00 & $1.11$e-02 & $99.6$ & $5.8$e-07 & $3.4$e-06 \\
        $120$ & $2394$ & $4851$ & $1.11$e+00 & $1.2$e-04 & $99.9$ & $3.65$e+00 & $1.57$e-02 & $99.6$ & $1.6$e-07 & $6.5$e-07 \\
        $140$ & $2653$ & $6843$ & $1.72$e+00 & $1.2$e-04 & $99.9$ & $4.57$e+00 & $2.25$e-02 & $99.5$ & $5.3$e-08 & $1.9$e-07 \\
        \hline
    \end{tabular}
    \endgroup
    \caption{{\em\small Timings (in seconds) for an experiment with $p=5$, grid size $h = 9/N$; the number of degrees of freedom and smooth quadrature sources can be determined (using $q = 2p$) from Equations~\eqref{eq:ndofs} and~\eqref{eq:nsrcs}, respectively, using the columns labeled $N_t$ and $N_b$. The quantities $T_\Gamma^\mathrm{fmm}$ and $T_\Gamma^\mathrm{corr}$ relate to evaluation of $\mathcal{V}[f]$ on $\Gamma$ at boundary integral collocation nodes (cf.~\eqref{eq:ellipt_augment_prob_bc}) while $T_\Omega^\mathrm{fmm}$ and $T_\Omega^\mathrm{corr}$ relate to evaluation of $u$ at all $\mathrm{ndofs}$ source function nodes in the domain; `fmm' superscripts refer to the application of the smooth quadrature rule and `corr' superscripts refer to the application of the sparse correction map. The time $T^\mathrm{bie}_\Gamma$ associated with solution of the boundary integral equations using GMRES, evaluation of the solution at all volumetric targets, and application of close corrections to the layer potential evaluations consistently required $T^\mathrm{bie}_\Gamma = 1.4$. The $\%$fmm figures are determined by the ratios of $T_\Gamma^\mathrm{fmm}$ to $T_\Gamma^\mathrm{corr}$ and of $T_\Omega^\mathrm{fmm}$ to $T_\Omega^\mathrm{corr}$, and a percentage of $99.9$ is displayed for any percentage above $99.9$.}} \label{tab:fmm_timing}
\end{table}

\section{Conclusions \& future directions}\label{sec:Conclusions}
We have presented a volume quadrature scheme that is used to solve inhomogeneous PDEs (either exterior or interior problems) on irregular (either simply- or multiply-connected) domains, and which is to our knowledge the first scheme for complex geometry that is of optimal asymptotic complexity in addition to being high-order accurate (including near to the boundary), and which, further, is capable of producing solutions to a wide variety of inhomogeneous elliptic PDEs, is efficient with respect to the number of degrees of freedom, is amenable to rapid generation of high-order accurate quadrature rules for singular and near-singular corrections, and is readily generalizable to non-PDE kernels and to three dimensional domains. It is worth mentioning that the present work solves inhomogeneous PDEs by direct evaluation of the true volume potential over $\Omega$, in contrast to the extension-based methods that constitute most prior work. In certain circumstances this distinction is relevant because the true volume potential is desired in and of itself and extension-based methods are no longer helpful; we mention for example ongoing work in volume potentials that arise in fractional PDEs and in the physical problem of quantifying fluid mixing, as well as the potential-theoretic solution, via the Lippmann-Schwinger equation, of the problem of wave scattering by media with a spatially-variable (and potentially discontinuous) index of refractivity.

A few direct extensions are currently being pursued, among them adaptivity in the bulk region and coupling to meshing technologies such as \texttt{TriWild} which could enable further adaptivity in the region near to the boundary. Use of the methodology with non-PDE kernels is currently being investigated, as is application to nonlinear PDEs. In the context of PDEs with moving geometries, a principal concern is the rapid generation of new quadratures for cells of the irregular section of the grid which naturally changes; the approach to efficient pre-computation of singular corrections we described should ameliorate this cost in that context. Finally, this work presents a clear path to tackling the 3D volume potential problem since all elements of the proposed methodology generalize immediately to that setting, with quadratures and orthogonal polynomial systems on (mapped curvilinear) triangles translating to (mapped curvilinear) tetrahedrons.

\section{Acknowledgements}
We acknowledge support from NSF under grant DMS-2012424, and the Mcubed program at the University of Michigan (UM). Authors also acknowledge the computational resources and services provided by UM's Advanced Research Computing. The authors are grateful to Zydrunas Gimbutas for helpful discussions.

\appendix
\begin{appendices}
\section{Proofs of supporting lemmas}
\begin{proof}[Proof of \Cref{lem:smooth_poincare}]
    Because $\mathcal{S}$ is star-shaped with respect to the origin the function $\rho(t\vv{r})$ is defined for all $\vv{r} \in \partial \mathcal{S}$ and for all $t \in [0, 1]$, and as a consequence the right-hand integral in~\eqref{eq:smooth_poincare} is well-defined.
    Denoting $\vv{r} = (x, y)$ and, for $\vv{r} \in \partial \mathcal{S}$, $\vv{\tau} = (dx/ds,dy/ds)$ (which is well-defined at almost every $\vv{r} \in \partial \mathcal{S}$ since $\partial\mathcal{S}$ is piecewise smooth), we first have
    \begin{align*}
            \oint_{\partial \mathcal{S}}\left(\int_0^1 t\,\rho(t\vv{r})\,\mathrm{d}t\right)  \vv{r}\times \vv{\tau}\,\mathrm{d}s
            =& \oint_{\partial \mathcal{S}} \left(\int_0^1 t\,\rho(t\vv{r})\,\mathrm{d}t\right)  \left(x \,\mathrm{d}y - y\,\mathrm{d}x \right).
    \end{align*}
    But then, since $\mathcal{S}$ is a region with a piecewise-smooth boundary and the function $\int_0^1 t\rho(t\vv{r})\,\mathrm{d} t$ is continuously differentiable with respect to each of the variables $x$ and $y$, we can apply Green's theorem~\cite[Thm.\ 10--43]{Apostol} to the vector field in the circulation integral to obtain
    \begin{align*}
            \oint_{\partial \mathcal{S}}\left(\int_0^1 t\,\rho(t\vv{r})\,\mathrm{d}t\right)  \vv{r}\times \vv{\tau}\,\mathrm{d}s
            =& \int_{\mathcal{S}} \left(\frac{\partial}{\partial x}\left(x\int_0^1 t\,\rho(t\vv{r})\,\mathrm{d}t\right) + \frac{\partial}{\partial y}\left(y\int_0^1 t\,\rho(t\vv{r})\,\mathrm{d}t\right)\right)\,\mathrm{dA}.
    \end{align*}
    It is then a straightforward computation to verify that
    \begin{align*}
            \oint_{\partial \mathcal{S}}\left(\int_0^1 t\,\rho(t\vv{r})\,\mathrm{d}t\right)  \vv{r}\times \vv{\tau}\,\mathrm{d}s =&
            \int_{\mathcal{S}}  \left(\int_0^1\left( \frac{\mathrm{d}}{\mathrm{d}t}t^2\,\rho(t\vv{r})\right)\,\mathrm{d}t \right)\,\mathrm{dA}
            = \int_{\mathcal{S}} \rho(\vv{r})\,\mathrm{dA},
    \end{align*}
    which completes the proof.
\end{proof}

\begin{proof}[Proof of \Cref{lem:poincare}]
    The proof proceeds by isolating a ball containing the location of the singularity, showing the contribution due to this region is negligible, and applying \Cref{lem:smooth_poincare} on a region excluding this ball. In order to proceed and utilize \Cref{lem:smooth_poincare} we first change to a coordinate system centered at $\vv{r}_*$ via the map $\vv{r} \mapsto \vv{r} + \vv{r}_*$, and write
    \begin{equation}\label{eq:S_rstar_translate}
        \int_{\mathcal{S}} K(\vv{r} - \vv{r}_*) \rho(\vv{r}) dA = \int_{\mathcal{S}_*} K(\vv{r}) \rho(\vv{r} + \vv{r}_*) dA,
\end{equation}
    where $\mathcal{S}_*$ is an $\vv{r}_*$-translate of $\mathcal{S}$ and, in particular, is star-shaped with respect to the origin.

    It is useful to introduce a smooth cut-off function $\rchi: [0, \infty) \to [0, 1]$, $\rchi$ non-increasing and satisfying both $\rchi(v) = 1$ in a neighborhood of $v = 0$ and $\rchi(v) = 0$ for $v \ge 1$. In view of~\eqref{eq:S_rstar_translate}, using this function we have the identity, for every $\varepsilon > 0$,
    \begin{equation}\label{eq:lemma_mollifier_split}
        \begin{split}
            \int_{\mathcal{S}} K(\vv{r})\rho(\vv{r} + \vv{r}_*)\,\mathrm{dA} = & \int_{\mathcal{S}_*} K(\vv{r})\left(1-\rchi(|\vv{r}|/\varepsilon)\right)\rho(\vv{r} + \vv{r}_*)\,\mathrm{dA}\\
            & + \int_{\mathcal{S}_*} K(\vv{r})\rchi(|\vv{r}|/\varepsilon)\rho(\vv{r} + \vv{r}_*)\,\mathrm{dA}.
        \end{split}
    \end{equation}
    It is easy to see that since for sufficiently small $\varepsilon > 0$ the quantity $\left(1-\rchi(|\vv{r}|/\varepsilon)\right)$ vanishes for $\vv{r}$ in an $\varepsilon$-dependent neighborhood of the origin, the integrand of the first integral on the right-hand side of~\eqref{eq:lemma_mollifier_split} is a smooth function of $\vv{r}$, and so by applying \Cref{lem:smooth_poincare} to that integral we have
\begin{equation}
\begin{split}
    \int_{\mathcal{S}_*} K(\vv{r})\rho(\vv{r} + \vv{r}_*)\,\mathrm{dA} = &\oint_{\partial \mathcal{S}_*}\left(\int_0^1 t\,K(t\vv{r})\left(1-\rchi(t|\vv{r}|/\varepsilon)\right)\rho(t\vv{r} + \vv{r}_*)\,\mathrm{d}t\right) \vv{r}\times \vv{\tau}\,\mathrm{d}s\\
    &+ \int_{\mathcal{S}_*} K(\vv{r})\rchi(|\vv{r}|/\varepsilon)\rho(\vv{r} + \vv{r}_*)\,\mathrm{dA}.\label{eq:lemma_mollifier_poincare}
\end{split}
\end{equation}
    Denoting by $I_\varepsilon(\vv{r})$ the integrand of the outer integral above,
    \[
        I_\varepsilon(\vv{r}) = (\vv{r} \times \vv{\tau}) \int_0^1 tK(t\vv{r}) \left(1-\rchi(t|\vv{r}|/\varepsilon)\right)\rho(t\vv{r} + \vv{r}_*)\,\mathrm{d}t,
    \]
    defining $\vv{\hat{r}}$ via $\vv{r} = |\vv{r}| \vv{\hat{r}}$, and estimating using the triangle inequality and the fact that $0 \le \rchi(v) \le 1$ we find that the bound
    \[
        \left|I_\varepsilon(\vv{r})\right| \le |\vv{r}| |\vv{\hat{r}} \times \vv{\tau}| \int_0^1 t\left|K(t\vv{r})\right|\rho(t\vv{r} + \vv{r}_*)\,\mathrm{d}t \le C |\vv{r}| \int_0^1 t\left|K(t\vv{r})\right|\,\mathrm{d}t
    \]
    holds for every $\varepsilon > 0$ and for all $\vv{r} \in \mathcal{S}_*$, where $C = C(\mathcal{S}_*, \rho) > 0$ is a constant. Since $K$ is weakly-singular, we further have that there exists a constant $\tilde{C} = \tilde{C}(\mathcal{S}_*,\rho, K) > 0$ such that, for all $\vv{r} \in \mathcal{S}_*$,
    \begin{equation}\label{eq:Ieps_bound}
        \left|I_\varepsilon(\vv{r})\right| \le \tilde{C} \left|\vv{r}\right|^{1-\alpha}\int_0^1 t^{1-\alpha}\,\mathrm{d}t,
    \end{equation}
    for some $\alpha < 2$, from which we conclude that $I_\varepsilon(\vv{r})$ is integrable over $\mathcal{S}_*$. Now, since firstly the upper bound for $I_\varepsilon(\vv{r})$ is integrable and since the integrand $I_\varepsilon(\vv{r})$ converges pointwise, and secondly the inner integral in $I_\varepsilon(\vv{r})$ also converges pointwise and is bounded above a function that is also integrable over the interval $[0, 1]$ (namely, by the integrand in the upper bound~\eqref{eq:Ieps_bound}), by two applications of the dominated convergence theorem in the $\varepsilon \to 0$ limit for the first integrand in the right-hand side of~\eqref{eq:lemma_mollifier_poincare} we find
    \begin{equation}
        \label{eq:eps_zero_poincare_main}
        \lim_{\varepsilon \to 0} \int_{\mathcal{S}_*} K(\vv{r})\left(1-\rchi(|\vv{r}|/\varepsilon)\right)\rho(\vv{r} + \vv{r}_*)\,\mathrm{dA} = \oint_{\partial \mathcal{S}_*}\left(\int_0^1 t\,K(t\vv{r}) \rho(t\vv{r} + \vv{r}_*))\,\mathrm{d}t\right)  \left(\vv{r}\times \vv{\tau}\right)\,\mathrm{d}s.
    \end{equation}
    But since $K$ is a weakly-singular kernel function, i.e.\ it satisfies for some constant $D > 0$ the bound $|K(\vv{r})| \le D |\vv{r}|^{-\alpha}$, $\alpha < 2$, we have that
    \begin{equation}\label{eq:integrable}
        \lim_{\varepsilon\rightarrow 0} \int_{B_\varepsilon(\vv{0})} |K(\vv{r})|\,\mathrm{dA} = 0,
    \end{equation}
    so that by taking the limit as $\varepsilon \to 0$ in the right-hand side of~\eqref{eq:lemma_mollifier_split} and using both~\eqref{eq:eps_zero_poincare_main} and~\eqref{eq:integrable}, we have
    \begin{equation}\label{eq:smooth_poincare_shifted}
        \int_{\mathcal{S}} K(\vv{r}-\vv{r}_*) \rho(\vv{r}) \,\mathrm{dA} = \oint_{\partial \mathcal{S}_*}\left(\int_0^1 t\,K(t\vv{r}) \rho(t\vv{r} + \vv{r}_*)\,\mathrm{d}t\right)  \vv{r}\times \vv{\tau}\,\mathrm{d}s.
    \end{equation}
    The desired result follows by a change of coordinates $\vv{r} \mapsto \vv{r} - \vv{r}_*$.
\end{proof}
\end{appendices}

\printbibliography
\end{document}